\documentclass[twoside,a4paper,12pt]{article}
\usepackage{amsmath}

\usepackage[usenames,dvipsnames]{color}

\def\Hassan{\color{red}}
\def\pier{\color{blue}}
\def\pcgg{\color{red}}
\let\Hassan\relax
\let\pier\relax
\let\pcgg\relax

\newcommand\testopari{\sc P. Colli, M. H. Farshbaf-Shaker, G. Gilardi and J. Sprekels}
\newcommand\testodispari{\sc Optimal control of a double obstacle Cahn--Hilliard inclusion}
\date{\today}

\begin{document}
\begin{center}{{\huge {\bf Optimal boundary control\\[2mm]
of a viscous Cahn--Hilliard system\\[2mm] 
with dynamic boundary condition\\[3mm] 
and double obstacle potentials}}}\\[6mm]

\vspace{9mm}
{\bf Pierluigi Colli\footnote{Dipartimento di Matematica  ``F. Casorati'',
Universit\`a di Pavia, Via Ferrata, 1, 27100 Pavia,  Italy,
e-mail: pierluigi.colli@unipv.it, gianni.gilardi@unipv.it},
{\pier M. Hassan Farshbaf-Shaker}\footnote{Weierstrass Institute for 
Applied Analysis and Stochastics, 
Mohrenstrasse 39, 10117 Berlin, Germany,
e-mail: Hassan.Farshbaf-Shaker@wias-berlin.de, juergen.sprekels@wias-berlin.de\\
\hspace*{5mm}$^3$Department of Mathematics, Humboldt-Universit\"at zu Berlin, Unter 
den Linden 6, 10099 Berlin, Germany\\[2mm]
{\bf Acknowledgements:} This paper was initiated during a visit of JS to the 
Universit\`a di Pavia. The kind hospitality and stimulating atmosphere of the
Universit\`a di Pavia are 
gratefully acknowledged. Some financial support 
comes from the MIUR-PRIN Grant 2010A2TFX2 ``Calculus of Variations''.},\\[0.2cm]
{\pier Gianni Gilardi}$^{1}$ {\pier {\rm and} J\"urgen Sprekels}$^{2\,\,3}$
}
%\footnote{Weierstrass Institute for 
%Applied Analysis and Stochastics,
%Mohrenstrasse 39, 10117 Berlin, Germany,
%e-mail: juergen.sprekels@wias-berlin.de\\[2mm]
\end{center}

\renewcommand\theequation{\mbox{\arabic{section}.\arabic{equation}}}
\newtheorem{theorem}{Theorem}[section]
\mathsurround .25mm
\newcommand{\calY}{{\mathcal Y}}
\newcommand{\dega}{{\Delta_\Gamma}}
\newcommand{\buuga}{{\bar u_\Gamma}}
\newcommand{\nf}{{\bf n}}
\newcommand{\ugesalph}{(u^\alpha,u_\Gamma^\alpha)}
\newcommand{\bugesalph}{(\bar{u}^\alpha,\bar{u}_\Gamma^\alpha)}
\newcommand{\ygesalph}{(y^\alpha,y_\Gamma^\alpha)}
\newcommand{\bygesalph}{(\bar{y}^\alpha,\bar{y}_\Gamma^\alpha)}
\newcommand{\uga}{{u_\Gamma}}
\newcommand{\uuga}{{(u,u_\Gamma)}}
\newcommand{\yga}{{y_\Gamma}}
\newcommand{\vga}{{v_\Gamma}}
\newcommand{\xiga}{{\xi_\Gamma}}
\newcommand{\rz}{{\rm I\!R}}
\newcommand{\nz}{{\rm I\!N}}
\newcommand{\ve}{{\varepsilon}}
\newcommand{\dx}{{{\rm d}x}}
\newcommand{\dt}{{{\rm d}t}}
\newcommand{\ds}{{{\rm d}s}}
\newcommand{\dgm}{{{\rm d}\Gamma}}
\newcommand{\essinf}{\mathop{\rm ess\, inf}}
\newcommand{\esssup}{\mathop{\rm ess\, sup}}
\newcommand{\oma}{{\Omega}}
\newcommand{{\tinto}}{{\int_0^T}}
\newcommand{{\xinto}}{{\int_\Omega}}
\newcommand{{\ginto}}{{\int_\Gamma}}
\newcommand{{\txinto}}{{\int_0^T\!\!\int_\Omega}}
\newcommand{{\tgamma}}{{\int_0^T\!\!\int_\Gamma}}
\newcommand{\txt}{{\int_0^t\!\!\int_\Omega}}
\newcommand{\txg}{{\int_0^t\!\!\int_\Gamma}}
\newcommand{\lzo}{{L^2(\Omega)}}
\newcommand{\heins}{{H^1(\Omega)}}
\newcommand{\lio}{{L^\infty(\Omega)}}
\newcommand{\glio}{{L^\infty(\Gamma)}}
\newcommand{\gheins}{{H^1(\Gamma)}}
\newcommand{\qlzo}{{L^2(Q)}}
\newcommand{\qlio}{{L^\infty(Q)}}
\newcommand{\glzsig}{{L^2(\Sigma)}}
\newcommand{\glisig}{{L^\infty(\Sigma)}}
\newcommand{\uad}{{\cal U}_{\rm ad}}
\newcommand{\qed}{\hfill\colorbox{black}{\hspace{-0.01cm}}}
%{\parbox{0.4cm}}}
\renewcommand{\min}{\mathop{\rm Min}}

% Das Layout wird mit Hilfe des geometry-Pakets erstellt.
% Eine einfache Aenderung der Textbreite und -hoehe bei Beibehaltung eines
% symmetrischen Blattaufbaus erhaelt man mit dem Kommando
% \geometry{twoside=false,textwidth=15cm,textheight=9in}

\title{Optimal boundary control of a viscous Cahn--Hilliard equation with dynamic boundary conditions and double obstacles}

\vspace{5mm}
\begin{abstract} 
\noindent
{\small
In this paper, we investigate optimal boundary control problems for Cahn--Hilliard variational inequalities with a dynamic boundary condition involving double obstacle potentials and the Laplace--Beltrami operator.  The cost functional is of standard tracking type, and box constraints for the controls are prescribed. We prove existence of optimal controls and derive first-order necessary conditions of 
optimality. The general strategy, which follows the lines of the recent 
approach {\pier by Colli, Farshbaf-Shaker, Sprekels (see Appl. Math. Optim., 2014)} to the (simpler) Allen--Cahn case, is the following: we use the results that were recently established by {\pier Colli, Gilardi, Sprekels 
in the preprint arXiv:1407.3916 [math.AP]} for the case of (differentiable) logarithmic potentials and perform a so-called ``deep quench limit''. Using compactness and monotonicity arguments, it is shown that this strategy leads to the desired first-order necessary optimality conditions for the case of (non-differentiable) double obstacle potentials.} 
\end{abstract}

%%%%%%%%%%%%%%%%%%%%%%%%%%%%%%%%%%%%%%%%%%%%%%%%%%%%%%%%%%%%%%%%%%%%%%
%                                                                    %
%                         Platz fuer neue Theoreme                   %
%                                                                    %
%%%%%%%%%%%%%%%%%%%%%%%%%%%%%%%%%%%%%%%%%%%%%%%%%%%%%%%%%%%%%%%%%%%%%%
\section{Introduction}
\noindent 
Let $\oma\subset\rz^N$, $2\le N\le 3$, denote some open{\pier , connected} and bounded
domain with smooth boundary $\Gamma$ and outward unit normal field $\nf$, and let $T>0$ be a
fixed final time. Putting $Q:=\oma\times (0,T)$, $\Sigma:=\Gamma\times (0,T)$, we introduce the function spaces

\begin{eqnarray}
\label{eq:1.1}
&&H:=L^2(\oma),\quad V:=H^1(\oma),\quad H_\Gamma:=L^2(\Gamma), \quad V_\Gamma:=H^1(\Gamma),\nonumber\\[2mm]
&&{\cal H}:=H\times H_\Gamma,\quad {\cal V}:=\left\{(y,y_\Gamma):\,y\in V\times V_\Gamma:
v_{{\pier |_\Gamma}}=v_\Gamma\right\},
\end{eqnarray}
which are Hilbert spaces when endowed with the topolgies induced by their respective natural inner products, {\pier denoted} by 
$\,(\,\cdot\,,\,\cdot\,)_E\,$ for $\,E\in \{H,H_\Gamma,V,V_\Gamma,{\cal H},{\cal V}\}$. 
In the following, we denote the norm in the generic Banach
space $E$ by $\|\,\cdot\,\|_E$, with the one exception that for convenience the norm of the space $H^N$ will 
also be denoted by $\|\,\cdot\,\|_H$. Moreover, {\pier let $\,E^*\,$ indicate the dual space of $\,E\,$ and let $\langle\,\cdot\,,\,\cdot\,\rangle_E\,$ always stand} for the duality pairing between elements of $\,E^*\,$ 
and elements of $\,E$.
It is understood that $H$ is embedded in $V^*$
in the usual way, namely, such that $\,\langle\,u\,,\,v\,\rangle_V=(\,u\,,\,v)_H\,$ for all $u\in H$ and $v\in V$;
we then obtain the Hilbert triplet $V\subset H\subset V^*$ with dense and compact 
embeddings. In the same way, we construct the Hilbert triplets $V_\Gamma\subset H_\Gamma
\subset V_\Gamma^*$ and ${\cal V}\subset {\cal H}\subset {\cal V}^*$, with dense and compact embeddings.

\vspace{2mm}
Throughout this paper, we generally assume:

\vspace{5mm}
\noindent
{\bf (A1)}\quad There are given constants $\beta_i\ge 0$, $1\le i\le 5$, which do not all vanish, as well as functions
\begin{eqnarray*}
&&z_Q\in L^2(Q),
\quad z_\Sigma\in L^2(\Sigma),\quad z_\Omega\in L^2(\oma),
\quad z_{\Gamma}\in L^2(\Gamma),  \quad\mbox{and}\\[1mm]
&&\widetilde{u}_{1_\Gamma}, \widetilde{u}_{2_\Gamma}\in L^\infty(\Sigma) 
\mbox{\,\, with \,\,}
\widetilde{u}_{1_\Gamma} \,\le \,\widetilde{u}_{2_\Gamma}
\mbox{\,\, a.\,e. on \,}\,\Sigma\,.
\end{eqnarray*}

\vspace{2mm}
We then introduce the tracking type cost functional 
\begin{eqnarray}
\label{eq:1.2}
&&{\cal J}((y,\yga),\uga):=\frac {\beta_1} 2 \,\|y-z_Q\|^2_{\qlzo}\,+\,\frac {\beta_2} 2\,
\|\yga-z_\Sigma\|^2_{L^2(\Sigma)} \nonumber\\[2mm]
&&\quad +\,\frac{\beta_3}2 \|y(T)-z_\Omega\|^2_{\lzo}\, +
\,\frac {\beta_4} 2 \,\|\yga(T)-z_{\Gamma}\|^2_{L^2(\Gamma)}
+\,\frac{\beta_5}2 \|u_\Gamma\|^2_{L^2(\Sigma)}\,,\quad
\end{eqnarray}
which is meaningful for, e.\,g., $(y,y_\Gamma)\in {\cal V}$ and $\uga\in H_\Gamma$,
and, for $\tau>0$, the viscous Cahn--Hilliard system with dynamic boundary conditions

\vspace{2mm}
\begin{equation}
\label{eq:1.3}
\partial_t y-\Delta w=0 \,\quad\mbox{in }\,Q,
\end{equation}
\begin{equation}
\label{eq:1.4}
w=\tau\,\partial_t y-\Delta y +\xi + f'_2(y) \,\quad\mbox{in }\,Q,
\end{equation}
\begin{equation}
\label{eq:1.5}
\quad y_{|_\Gamma}=\yga,\quad  \partial_\nf y \,+\, \partial_t y_{\Gamma}-\dega \yga +\xiga +  g_{2}'(\yga)=\uga\,,\quad\partial_\nf w=0\,,\quad\mbox{on }\,\Sigma,
\end{equation}
\begin{equation}
\label{eq:1.6}
\xi\in\partial I_{[-1,1]}(y)\quad\mbox{a.\,e. in }\,Q,\quad \xiga\in\partial I_{[-1,1]}(y_\Gamma)\quad\mbox{{\pier a.\,e.} on }\,\Sigma,
\end{equation}
\begin{equation}
\label{eq:1.7}
y(\cdot,0)=y_0 \,\quad\mbox{a.\,e. in }\,\oma,\,\quad\,y_\Gamma(\cdot,0)=y_{0_\Gamma}\,\quad
\mbox{a.\,e. on }\,\Gamma\,.
\end{equation}

\vspace{2mm}
Moreover, let $M_0>0$ denote some given constant, and let 
\begin{eqnarray}
\label{eq:1.8}
&&\uad := \left\{u_\Gamma\in H^1(0,T;H_\Gamma)\cap L^\infty(\Sigma)\,:\,
\|\partial_t \uga\|_{L^2(\Sigma)}\,\le\,M_0\right.\,,\nonumber\\
&&\left.\qquad\qquad \widetilde{u}_{1_\Gamma}\le u_\Gamma\le
\widetilde{u}_{2_\Gamma} \quad {\mbox{a.\,e. in }\,\Sigma} \ \right\}, 
\end{eqnarray} 
be the set of admissible controls which is assumed nonempty throughout this paper.
Our overall boundary control problem reads as follows:
\vspace{2mm} 
\begin{eqnarray*}
({\mathcal{P}}_{0})\quad&&\mbox{Minimize }\,\, {\cal J}((y,\yga),\uga)
\quad\mbox{subject to the state constraints (\ref{eq:1.3})--(\ref{eq:1.7})} \\
&&\mbox{{\pier and to the} control constraint }\,\uga\in \uad.
\end{eqnarray*}

\vspace{2mm}
In (\ref{eq:1.7}), $y_0$ and $y_{0_\Gamma}$ are given initial data with $y_{0|_{\Gamma}}=y_{0_\Gamma}$, where the trace $\,y_{|_\Gamma}\,$ (if it exists) of a function $y$ on $\Gamma$ will throughout be denoted by $y_\Gamma$ without further comment. Moreover, 
in the following $\partial_\nf$, $\nabla_\Gamma$ and $\dega$ will always stay for
the outward normal derivative, the tangential gradient, and the Laplace--Beltrami 
operator, respectively, on $\Gamma$; in addition, $f_{2}\,,\,g_2$ are given smooth nonlinearities, while $\uga$ is a boundary control. Since we will confine ourselves to the   viscous case $\tau>0$, we will henceforth assume without loss of generality
that $\tau=1$. 

\vspace{2mm}
The system (\ref{eq:1.3})--(\ref{eq:1.7}) is an initial-boundary 
value problem with nonlinear dynamic boundary condition for a Cahn--Hilliard differential inclusion, which (cf. Proposition 2.2 below) under appropriate conditions on the data  admits for every $\uga\in\uad$ a solution quintuple $(y,\yga,w,\xi,\xiga)$, where the solution
components $(y,\yga,\xiga)$ are uniquely determined. Hence, the {\em control-to-state 
operator} ${\cal S}_0: \uga\mapsto {\cal S}_0(u_\Gamma):=(y,\yga)$  is well defined on $\uad$, and the control problem $({\cal P}_0)$ is equivalent to minimizing the
reduced cost functional
\begin{equation}
\label{eq:1.9}
{\cal J}_{\rm red} (\uga) \,:=\, {\cal J}({\cal S}_0(\uga),\uga)
\end{equation}
over $\uad$.

\vspace{2mm}
In the physical interpretation, the unknown $y$ usually stands for the (conserved) order parameter of an
isothermal phase transition, typically a rescaled fraction of one of the involved phases. 
{\pier In such a situation,} it is physically meaningful to require $y$ to attain values in the interval $[-1,1]$ on
both $\oma$ and $\Gamma$. A standard technique to meet this requirement is to use the indicator function of the interval $[-1,1]$,
$$I_{[-1,1]}(y)=\left \{
\begin{array}{ll}
0&\mbox{if }\,y\in [-1,1]\\
+\infty&\mbox{otherwise}
\end{array} \right. ,
$$
 so that the non-diffusive parts of the local specific bulk  and surface free energies, 
$F_{\rm bulk}:=I_{[-1,1]}+f_2$ and  $F_{\rm surface}:=
I_{[-1,1]}+g_{2}$, are of {\em double obstacle type}, 
and the subdifferential 
$\partial I_{[-1,1]}$, defined~by 
$$
 \eta \in \partial I_{[-1,1]}( v)  \quad \hbox{ if and only if }
\quad
\eta \ \left\{
\begin{array}{ll}
\displaystyle
\leq \, 0 \   &\hbox{if } \ v=-1  
\\[0.1cm]
= \,0\   &\hbox{if } \ -1 < v <1  
\\[0.1cm]
\geq \, 0 \  &\hbox{if } \  v =1  
\\[0.1cm]
\end{array}
\right. ,
$$
is employed in place of the usual derivative.
Concerning the selections $\,\xi$, $\xiga$\, in 
(\ref{eq:1.6}), one has to keep in mind that $\xi$ 
may be not regular enough as to single  out its trace on 
the boundary $\Gamma$, and if the trace 
$ \xi_{|_\Gamma}$ exists, it may differ from $\xi_{\Gamma}$, in general.

\vspace{2mm}
\noindent The optimization problem $({\mathcal{P}}_{0})$ belongs to the problem class of so-called MPECs (Mathematical Programs with Equilibrium Constraints). It is a well-known fact that the differential inclusion conditions encoded in (\ref{eq:1.3})--(\ref{eq:1.6}), which  occur as constraints in $({\mathcal{P}}_{0})$, violate  all 
of the known classical nonlinear programming constraint qualifications. Hence, the existence of Lagrange multipliers cannot be inferred from standard theory, and the derivation of first-order necessary condition becomes very 
difficult.

\vspace{2mm}
While numerous papers deal with the well-posedness and asymptotic behavior of Cahn--Hilliard
system (cf., e.\,g., the references given in {\pier \cite{GiMiSchi, GiMiSchi2, CGS1}}),
there are comparatively few investigations of associated optimal control problems. Usually,
these papers treat the non-viscous case $\tau=0$ and are restricted to differentiable 
free energies and to the case of distributed controls, with the no-flux condition
$(\partial_\nf y)_{{\pier |_\Gamma}}=0$ assumed in place of the more difficult dynamic boundary condition (\ref{eq:1.5}).  In this connection, we refer to \cite{WaNa} and \cite{HW1}, where the latter paper also deals with the case of double obstacle potentials. 

\vspace{2mm}
Quite recently, also convective Cahn--Hilliard systems have been investigated from
the viewpoint of optimal control. In this connection, we refer to \cite{ZL1} and
\cite{ZL2}, where the latter paper deals with the two-dimensional case. The three-dimensional case with a nonlocal free energy was studied in \cite{RoSp}. There also exist 
contributions dealing with the more general and difficult Cahn--Hilliard/Navier--Stokes 
systems, cf.~\cite{HK} and \cite{HW2}. Finally, we mention the papers \cite{CGPS} and 
\cite{CGScon}, in which control problems for a generalized Cahn--Hilliard system
introduced in \cite{Podio} were investigated.

\vspace{2mm}
The only existing contribution to the optimal control of viscous or non-viscous 
Cahn--Hilliard systems with dynamic boundary conditions of the form (\ref{eq:1.5}) 
seems to be the recent paper \cite{CGS2} in which three of the present authors investigated 
the case of differentiable bulk and surface free energies that may have singular 
derivatives. A typical case to which the analysis in \cite{CGS2} applies is given by
the logarithmic form
\begin{eqnarray}
\label{eq:1.10}
F_{\rm log}(y)&\!\!=\!\!&h(y)+f_2(y)\,,\quad\mbox{where }\nonumber\\
h(y)&\!\!=\!\!& \widehat{c}\,((1+y)\,\ln(1+y)+(1-y)\,\ln(1-y)), \quad -1<y<1,
\end{eqnarray}
with some fixed constant $\widehat{c}>0$. Note that in this case the inclusions
(\ref{eq:1.6}) have to be replaced by the equations $\,\xi=h'(y)\,$ and $\,\xi_\Gamma
=h'(\yga)$, respectively.

\vspace{2mm} In this paper, we aim to employ the results established in \cite{CGS2}
to treat the non-differentiable double obstacle case when $\xi,\xi_\Gamma$ satisfy the 
inclusions (\ref{eq:1.6}). Our approach is guided by the strategy used by three 
of the present authors in their recent paper \cite{CFS} for a corresponding optimal
control problem for the simpler Allen--Cahn equation: in \cite{CFS}, necessary
optimality conditions for the double obstacle case could be established by 
performing a so-called ``deep quench limit'' in a family of optimal control problems
with differentiable nonlinearities of a form that had been previously treated in \cite{CS}
and for which the {\pier corresponding systems} had been analyzed in \cite{CC}.  

\vspace{2mm}
The general idea is briefly explained as follows: we replace the inclusions (\ref{eq:1.6})
by
\begin{equation}
\label{eq:1.11}
\xi=\varphi(\alpha)\,h'(y), \quad \xi_\Gamma=\psi(\alpha)\,h'(y),
\end{equation}
where $h$ is defined in (\ref{eq:1.10}), and where 
 $\varphi,\psi$ are continuous and positive functions on $(0,1]$ that satisfy 
\begin{equation}
\label{eq:1.12}
\lim_{\alpha\searrow 0}\,\varphi(\alpha)=\lim_{\alpha\searrow 0}\,\psi(\alpha)
=0, \quad \varphi(\alpha)\le C_{\varphi \psi}\,\psi(\alpha)\,\,\,\forall\,
\alpha>0,\,\,\mbox{with some }\, C_{\varphi\psi}>0.
\end{equation}
We remark that we could simply choose $\,\varphi(\alpha)
=\psi(\alpha)=\alpha^p\,$ for some $\,p>0$; however, there might be situations (e.\,g., in the
numerical approximation) in which it is advantageous to let $\,\varphi\,$ and $\,\psi\,$ have a different behavior
as $\,\alpha\searrow 0$.

\vspace{2mm} 
Now observe that $h'(y)=\ln\left(\frac{1+y}{1-y}\right)$ \,and\, $h''(y)=\frac 2 {1-y^2}>0$\, for 
$y\in (-1,1)$. Hence, in particular, we have
\begin{eqnarray}
\label{eq:1.13}
&&\lim_{\alpha\searrow 0}\,\varphi(\alpha)\,h'(y)=0 \quad\mbox{for }\, -1<y<1,\nonumber\\[2mm]
&&\lim_{\alpha\searrow 0}\Bigl(\varphi(\alpha)\,\lim_{y\searrow -1}h'(y)\Bigr)\,=\,-\infty,
\quad \lim_{\alpha\searrow 0}\Bigl(\varphi(\alpha)\,\lim_{y\nearrow +1}h'(y)\Bigr)\,=\,+\infty\,.
\end{eqnarray}
Since similar relations hold if $\varphi$ is replaced by $\psi$, we may regard the graphs of the functions 
$\,\varphi(\alpha)\,h'\,$ and $\,\psi(\alpha)\,h'\,$ as approximations to the graph of the subdifferential
$\partial I_{[-1,1]}$. 

\vspace{2mm}
Now, for any $\alpha>0$ the optimal control problem (later to be denoted by $({\cal P}_\alpha)$), which results if in $({\cal P}_0)$ the relation (\ref{eq:1.6}) is replaced by 
(\ref{eq:1.11}), is of the type for which
in \cite{CGS2} the existence of optimal controls $u_\Gamma^\alpha\in\uad$ as well as first-order necessary optimality conditions have been derived. Proving a priori estimates (uniform in $\alpha>0$), and 
employing compactness and monotonicity arguments, we will be able to show the following existence and approximation result: whenever $\,\{u_\Gamma^{\alpha_n}\}\subset\uad$ is a sequence of optimal controls for $({\cal P}_{\alpha_n})$, where $\alpha_n\searrow 0$ as $n\to\infty$, then there exist
a subsequence of $\{\alpha_n\}$, which is again indexed by $n$, and an optimal control 
$\bar u_\Gamma\in\uad$ of
$({\cal P}_0)$ such that

\begin{equation}
\label{eq:1.14}
u_\Gamma^{\alpha_n}\to\bar u_\Gamma \quad\mbox{weakly-star in ${\cal X}$ as }\,
n\to\infty,
\end{equation}
where here and in the following
\begin{equation}
\label{eq:1.15} 
{\cal X}:=H^1(0,T;H_\Gamma)\cap L^\infty(\Sigma)
\end{equation}
will always denote the control space.
In other words, optimal controls for $({\cal P}_\alpha)$ are for small $\alpha>0$ likely to be `close' to 
optimal controls for $({\cal P}_0)$. It is natural to ask if the reverse holds, i.\,e., whether every optimal control for
 $({\cal P}_0)$ can be approximated by a sequence $\,\{u_\Gamma^{\alpha_n}\}\,$ of optimal controls
for $({\cal P}_{\alpha_n})${\pier , for some} sequence $\alpha_n\searrow 0$. 

\vspace{2mm}
Unfortunately, we will not be able to prove such a `global' result that applies to all optimal controls for
$({\cal P}_0)$. However,  a `local' result can be established. To this end, let $\buuga\in\uad$ be any optimal control
for $({\cal P}_0)$. We introduce the `adapted' cost functional
\begin{equation}
\label{eq:1.16}
\widetilde{\cal J}((y,\yga),\uga) \,:=\,{\cal J}((y,\yga),\uga)\,+\,\frac 1 2
\|u_\Gamma-\bar u_\Gamma\|^2_{L^2(\Sigma)}
\end{equation}
and consider for every $\alpha\in (0,1]$ the {\em adapted control problem} of minimizing $\,\widetilde{\cal J}\,$ subject to $\uga\in\uad$ and to the constraint that $(y,\yga)$ solves the approximating system (\ref{eq:1.3})--(\ref{eq:1.5}), (\ref{eq:1.7}), 
(\ref{eq:1.11}). It will then turn out that the following is true: 

\vspace{2mm}
(i) \,There are some sequence $\,\alpha_n\searrow 0\,$ and minimizers 
$\,{\bar u_\Gamma^{\alpha_n}}\in\uad$ of the adapted control problem 
associated with $\alpha_n$, $n\in\nz$,
such that
\begin{equation}
\label{eq:1.17}
{\bar u_\Gamma^{\alpha_n}}\to\buuga\quad\mbox{strongly in $\glzsig$
as }\,n\to \infty.
\end{equation}
(ii) It is possible to pass to the limit as $\alpha\searrow 0$ in the first-order necessary
optimality conditions corresponding to the adapted control problems associated with $\alpha\in (0,1]$ in order to derive first-order necessary optimality conditions for problem $({\cal P}_0)$.

\vspace{2mm}
The paper is organized as follows: in Section~\ref{state}, we give a precise statement of the problem
under investigation, and we derive some results concerning the state system 
(\ref{eq:1.3})--(\ref{eq:1.7}) and 
its $\alpha$-approximation which is obtained if in $({\cal P}_0)$ the relations
 (\ref{eq:1.6}) are replaced by the relations (\ref{eq:1.11}).
In Section~\ref{existence}, we then prove the existence of optimal controls and the approximation result formulated above in
(i). The final Section~\ref{optimality} is devoted to the derivation of the first-order necessary 
optimality conditions, where the  strategy outlined in (ii) is employed. 

\vspace{2mm}
During the course of this analysis, we will make 
repeated use of the elementary Young's inequality
$$
a\,b\,\le\,\gamma |a|^2\,+\,\frac 1{4\gamma}\,|b|^2\quad\forall\,a,b\in\rz \quad\forall\,\gamma>0,
$$
and we will use the following notation: {\pier for functions $v\in V^*$ and 
$w\in L^1(0,T;V^*)$ we define their generalized mean values as}
\begin{equation}
\label{eq:1.18}
v^\oma:=\frac 1 {|\oma|}\,\langle \,v\,,\,1\,\rangle_V, \,\mbox{ and }\,
w^\oma(t):= (w(t))^\oma \,\mbox{ for a.\,e. }\,t\in (0,T).
\end{equation}
Clearly, (\ref{eq:1.18}) gives the usual mean values when {\pier elements of $H$ or
of $L^1(0,T;H)$, respectively, are involved}. We also recall Poincar\'{e}'s inequality
\begin{equation}
\label{eq:1.19}
\|v\|_V\,\le\,C_P\left(\|\nabla v\|_H\,+\,\left|v^\oma\right|\right)\quad\,\forall \,v\in V,
\end{equation} 
with a constant $C_P>0$ that only depends on $\oma$.
\section{General assumptions and state equations}\label{state}
\setcounter{equation}{0}
In this section, we formulate the general assumptions of the paper, and we state some preparatory results for the state system (\ref{eq:1.3})--{\pier(\ref{eq:1.7})} and its $\alpha$-approximations. 

\vspace{2mm}
We make the following general assumptions:
\vspace{2mm}

{\bf (A2)} \quad\,$f_2,g_2\in C^3([-1,1])$.\\[2mm]
{\bf (A3)} \quad\,{\pcgg $y_0\in H^2(\oma)$, $y_{0_\Gamma}:= {y_{0}}_{{\pier |_\Gamma}}\in H^2(\Gamma)$, 
% $\partial_\nf {y_{0}}_{{\pier |_\Gamma}}=0$, 
and we have} 
\begin{equation}
\label{eq:2.1}
-1<y_0(x)<1 \quad\forall\,x\in \overline{\oma} \,.
\end{equation}
{\bf (A4)} \quad\,There exist $\xi_0\in H$ and $\xi_{\Gamma,0}\in H_\Gamma$ such that
\begin{equation}
\label{eq:2.2}
\xi_0\in I_{[-1,1]}(y_0) \,\mbox{ a.\,e. in $\,\oma$}, \quad 
\xi_{\Gamma,0}\in I_{[-1,1]}(y_{0_\Gamma}) \,\mbox{ a.\,e. on $\,\Gamma$}.  
\end{equation}
 
\vspace{5mm} 
Now observe that the set $\uad$ is a bounded subset of ${\cal X}$. Hence, there exists
a bounded open ball in ${\cal X}$ that contains $\uad$. For later use it is convenient to fix such a ball once
and for all, noting that any other such ball could be used instead. In this sense, the following assumption
is rather a denotation:

\vspace{5mm}
{\bf (A5)}\quad${\cal U}$ is a nonempty open and bounded subset of ${\cal X}$ containing $\uad$, and the constant $\,R>0\,$ satisfies
\begin{equation}
\label{eq:2.3}
\|\uga\|_{H^1(0,T;H_\Gamma)}\,+\,\|u_\Gamma\|_\glisig\,\le\,R \quad\,\forall\,u_\Gamma\in {\cal U}.
\end{equation}  

\vspace{2mm}
Next, we introduce our notion of {\pier solution} to the problem (\ref{eq:1.3})--(\ref{eq:1.7}) in the abstract setting introduced above.

\vspace{5mm}
{\bf Definition 2.1:}\quad\,{\em A quintuple} $(y,\yga,w,\xi,\xi_\Gamma)$ {\em such that}
\begin{eqnarray}
\label{eq:2.4}
&&y\in H^1(0,T;H)\cap L^\infty(0,T;V)\cap L^2(0,T;H^2(\oma)),\\[1mm]
\label{eq:2.5}
&&\yga\in H^1(0,T;H_\Gamma)\cap L^\infty(0,T;V_\Gamma)\cap L^2(0,T;H^2(\Gamma)),\\[1mm]
\label{eq:2.6}
&&y \in [-1,1] \quad\mbox{\em a.\,e. in }\,Q, \quad \yga\in [-1,1]\quad\mbox{\em a.\,e. on }\,\Sigma,\\[1mm]
\label{eq:2.7}
&&\xi\in L^2(0,T;H) \,\,\mbox{\em and }\,\xi\in\partial I_{[-1,1]}(y) \,\,\,\mbox{\em a.\,e. in \,$Q$, }\,\\[1mm]
\label{eq:2.8} 
&&\xiga\in L^2(0,T;{\Hassan H_\Gamma}) \,\,\mbox{\em and }\,\xiga\in\partial I_{[-1,1]}(y_\Gamma)\,\,\,
\mbox{\em a.\,e. on $\Sigma$},\\[1mm]
\label{eq:2.9} 
&&w\in L^2(0,T;V),
\end{eqnarray}
{\pier {\em as well as} $\,\yga=y_{{\pier |_\Gamma}}$, ${y(0)}=y_0$, ${y_\Gamma(0)}=y_{0_\Gamma}$,} {\em is called a solution to} (\ref{eq:1.3})--(\ref{eq:1.7}) {\em if and only if it satisfies for almost every $t\in (0,T)$ the variational equations}
\begin{align}
\label{eq:2.10}
&\xinto \partial_ty(t)\,v\,\dx +\xinto\nabla w(t)\cdot\nabla v\,\dx =0 \quad\mbox{\em for
every $v\in V$},\\[2mm]
\label{eq:2.11}
&\xinto w(t)\,v\,\dx = \xinto \partial_t y(t)\,v\,\dx 
 +\xinto \nabla y(t)\cdot\nabla v\,\dx +\xinto(\xi(t)+f_2'(y(t)))\,v\,\dx 
\nonumber\\[2mm]
&\quad +\ginto\partial_t \yga(t)\,\vga\, \dgm+\ginto\nabla_\Gamma y_\Gamma(t)\cdot\nabla_\Gamma v_\Gamma\,\dgm\nonumber+\ginto(\xi_\Gamma(t)+g_2'(y_\Gamma (t))-\uga(t))\,{v_\Gamma}\,\dgm \nonumber\\[2mm]
&\quad\mbox{\em for every $\,(v,\vga)\in {\mathcal V}$.}
\end{align}

It is worth noting that (recall the notation (\ref{eq:1.18})) 
\begin{eqnarray}
\label{eq:2.12}
&&(\partial_t y(t))^\oma=0 \quad\mbox{for a.\,e. $t\in (0,T)$, and\, 
$y(t)^\oma=m_0$ \,for every $t\in [0,T]$,}\nonumber\\[2mm]
&&\mbox{where \,$m_0=(y_0)^\oma\,$ is the mean value of $y_0$,} 
\end{eqnarray} 
as usual for the Cahn--Hilliard equation. Notice that {\bf (A3)} implies $\,-1<m_0<1\,$ so that
$h'(m_0)$ is finite.

\vspace{5mm}
The following existence and uniqueness result follows from \cite[{\pier Theorems~2.2 and~2.4}]{CGS1}. {\pcgg Let us stress that the assumption (2.37) explicitely required in the statement of \cite[{\pier Thm.~2.4}]{CGS1} contains the condition $\partial_\nf {y_{0}}_{{\pier |_\Gamma}}=0$
which is completely useless (actually, it is never employed in the proof, as the reader can check).}

\vspace{5mm}
{\bf  Proposition~2.2:}\quad\,{\em Assume that {\bf (A2)}--{\bf (A4)} are fulfilled. Then there exists for any $\uga\in {\cal X}$ a  quintuple $(y,\yga,w,\xi,\xi_\Gamma)$ solving problem {\rm (\ref{eq:1.3})--(\ref{eq:1.7})} in the sense of Definition 2.1. For any such solution, we have the additional regularity
properties 
\begin{eqnarray*}
&&y\in W^{1,\infty}(0,T;H)\cap H^1(0,T;V)\cap L^\infty(0,T;H^2(\oma)),\\[1mm]
&&\yga\in W^{1,\infty}(0,T;H_\Gamma)\cap H^1(0,T;V_\Gamma)\cap L^\infty(0,T;H^2(\Gamma)),\\[1mm]
&&\xiga\in L^\infty(0,T;H_\Gamma).
\end{eqnarray*} 
Moreover, any two solution quintuples have the same components
$y, \yga, \xiga$ (while the components $w,\xi$ may not be uniquely determined).}

\vspace{5mm}
As in the Introduction, we denote the control-to-state operator, which assigns
to every \,$\uga\in {\mathcal X}\,$ the (uniquely determined) first two components $(y,\yga)$ of the associated solution quintuple, by $\,{\cal S}_0$.

\vspace{5mm}
We now turn our attention to the approximating state equations. As announced in the Introduction,  
we choose a special approximation of (\ref{eq:1.3})--(\ref{eq:1.7}); namely,
for $\alpha\in (0,1]$ we consider the system
 
\vspace{2mm} 
\begin{equation}
\label{eq:2.13}
\partial_t y^\alpha-\Delta w^\alpha = 0\,\quad\mbox{a.\,e. in $Q$}, 
\end{equation}
\begin{equation}
\label{eq:2.14}
w^\alpha={\partial_t y}^\alpha-\Delta y^\alpha\,+\,\varphi(\alpha)\,h'(y^\alpha)\,+\, f_{2}'(y^\alpha)\,\quad\mbox{a.\,e. in 
$Q$},
\end{equation}
\begin{align}
\label{eq:2.15}
\quad y^\alpha_{|_\Gamma}=y_\Gamma^\alpha,\quad \partial_\nf y^\alpha {}+{}
\partial_t y_{\Gamma}^\alpha-\dega y_\Gamma^\alpha  \,+\,\psi(\alpha)\, h'(y_\Gamma^\alpha) \,+\,  g_{2}'(y_\Gamma^\alpha)=\uga \,,\nonumber 
\\ {\pier \partial_\nf w^\alpha=0\,\quad\mbox{a.\,e. on $\Sigma$},}
\end{align}
\begin{equation}
\label{eq:2.16}
y^\alpha(\cdot,0)= y_0 \,\quad\mbox{a.\,e. in }\,\oma,\,\quad\,y_\Gamma^\alpha(\cdot,0)= y_{{0_\Gamma}} \,\quad
\mbox{a.\,e. on }\,\Gamma\,,
\end{equation}

\vspace{2mm}
where  $h$ is defined in (\ref{eq:1.10}) and $\varphi,\psi$ are positive and continuous functions on $(0,1]$ that
satisfy (\ref{eq:1.12}). Observe that as in (\ref{eq:2.10}), (\ref{eq:2.11}) the notion of a solution to 
(\ref{eq:2.13})--(\ref{eq:2.16}) has to be 
understood in the sense that for almost every $t\in (0,T)$ the following variational equations are satisfied:
\begin{align}
\label{eq:2.17}
&\xinto \partial_t y^\alpha(t)\,v\,\dx +\xinto\nabla w^\alpha(t)\cdot\nabla v\,\dx =0 \quad
{\pier \hbox{for every } \,} v\in V,\\[3mm]
\label{eq:2.18}
&\xinto\!\! w^\alpha(t)v\,\dx = \xinto\!\! \partial_t y^\alpha(t)v\,\dx 
 +\xinto\! \nabla y^\alpha(t)\cdot\nabla v\,\dx +\xinto(\varphi(\alpha)h'(y^\alpha(t))+f_2'(y^\alpha(t)))v\,\dx 
\nonumber\\[2mm]
&+\ginto\!\partial_t y_\Gamma^\alpha(t)\vga\, \dgm+\!\ginto\!\nabla_\Gamma y_\Gamma^\alpha(t)\cdot\nabla_\Gamma v_\Gamma\,\dgm\nonumber+\!\ginto\!(\psi(\alpha)h'(y_\Gamma^\alpha(t))+g_2'(y_\Gamma^\alpha (t))-\uga(t))\,{v_\Gamma}\,\dgm \nonumber\\[2mm]
&\quad{\pier \hbox{for every } \,}  (v,\vga)\in {\mathcal V}.
\end{align}

Since the functions $f^\alpha(y):=\varphi(\alpha)\,h(y)+f_2(y)$ and $f_\Gamma^\alpha(y):=
\psi(\alpha)\,h(y)+g_2(y)$ fulfill on $(-1,1)$ the conditions (2.3)--(2.7) in \cite{CGS2}, we can infer from
\cite[Thm.~2.1]{CGS2} that
the system (\ref{eq:2.13})--(\ref{eq:2.16}) admits for every $\uga\in {\cal U}$ a unique solution triple
$(y^\alpha, y_\Gamma^\alpha, w^\alpha)$  having the following properties: 
\begin{eqnarray}
\label{eq:2.19}
&&y^\alpha\in W^{1,\infty}(0,T;H)\cap H^1(0,T;V)\cap L^\infty(0,T;H^2(\oma)),\\[2mm]
\label{eq:2.20}
&&y_\Gamma^\alpha \in W^{1,\infty}(0,T;H_\Gamma)\cap H^1(0,T;V_\Gamma)\cap L^\infty(0,T;H^2(\Gamma)),\\[2mm]
\label{eq:2.21}
&&w^\alpha\in L^\infty(0,T;H^2(\oma)),\\[2mm]
\label{eq:2.22}
&&r_-^\alpha \,\le\,y^\alpha\,\le\,r_+^\alpha \,\,\mbox{ a.\,e in $Q$}, \qquad  
r_-^\alpha \,\le\,y^\alpha_\Gamma\,\le\,r_+^\alpha \,\,\mbox{ a.\,e on $\Sigma$},
\end{eqnarray}
with suitable constants $r_-^\alpha\,,\,r_+^\alpha \in (-1,1)$
that only depend on $\oma$, $T$, $y_0$, $y_{0_\Gamma}$, $f_2$, $g_2$, $\alpha$, and the constant $R>0$
introduced in {\bf (A5)}. In particular, the control-to-state mapping for the system
(\ref{eq:2.13})--(\ref{eq:2.16}), 
$\,{\cal S}_\alpha: \uga\mapsto {\cal S}_\alpha(\uga):=(y^\alpha,y_\Gamma^\alpha)$, for $\uga\in 
{\cal X}$, is well defined. Observe that the separation property (\ref{eq:2.22}) cannot be expected to
hold uniformly in $\alpha\in (0,1]$, in general; indeed, it cannot be excluded that there exists
some sequence $\{\alpha_n\}\subset (0,1]$ with $\alpha_n\searrow 0$ such that $r_-^{\alpha_n}\searrow -1$
and/or $r^{\alpha_n}_+ \nearrow +1$ as $n\to\infty$. 

\vspace{2mm}
We now aim to derive some
a priori estimates for $(y^\alpha,y_\Gamma^\alpha)$ which are independent of $\alpha$.
Prior to this, we recall a functional analytic framework which is customary in the context of
Cahn--Hilliard systems. We define
\begin{equation}
\label{eq:2.23}
\mbox{dom}\,{\cal N}:=\left\{v_*\in V^*: v_*^\oma=0\right\} \,\mbox{ and }\,
{\cal N}:\mbox{dom}\,{\cal N}\,\to\left\{v\in V: v^\oma=0\right\}
\end{equation} 
by setting for $v_*\in \mbox{dom}\,{\cal N}$
\begin{equation}
\label{eq:2.24}
{\cal N}v_*\in V, \quad ({\cal N}v_*)^\oma=0, \quad\mbox{and }\,\xinto\nabla{\cal N}v_*
\cdot\nabla z\,\dx=\langle v_*,z\rangle_V \,\quad\forall\,z\in V,
\end{equation}
that is, ${\cal N}v_*$ is the (unique) solution to the generalized Neumann problem 
$\,-\Delta v=v_*\,$ in $\oma$, $\partial_\nf v=0$ on $\Gamma$, {\pier that} satisfies $v^\oma=0$. 
Since $\,\oma\,$ is a bounded {\pier connected} domain with smooth boundary, it turns out that (\ref{eq:2.24})
yields a well-defined isomorphism that also {\pier fulfills}, for all $s\ge 0$,
\begin{align}
\label{eq:2.25}
{\cal N}v_*\in H^{s+2}(\oma) \,\,\mbox{ and }\,\, \|{\cal N}v_*\|_{H^{s+2}(\oma)}\,\le\,C_s\,
\|v_*\|_{H^s(\oma)} \nonumber \\[0.2cm] {\pier \mbox{for all }\, v_*\in H^s(\oma)\cap \mbox{dom}\,{\cal N},} 
\end{align}
where the constant $C_s>0$ depends only on $\oma$ and $s$. 
Moreover, if we define the mapping
$\,\|\,\cdot\,\|_*: V^*\to [0,+\infty)\,$ through the formula
\begin{equation}
\label{eq:2.26}
\|v_*\|_*^2:= \|\nabla{\cal N} (v_*-v_*^\oma)\|_H^2 \,+\,\left|v_*^\oma\right|^2 \quad
\forall\, v_*\in V^*, 
\end{equation} 
then it is straightforward to prove that $\|\,\cdot\,\|_*$ defines a norm on $V^*$ which turns
out to be equivalent to the usual norm of $V^*$. We thus have, with a constant $C_*>0$ that
depends only on $\oma$,
\begin{equation}
\label{eq:2.27}
\left|\langle v_*,v \rangle_V\right|\,\le\,C_*\,\|v_*\|_*\,\|v\|_V\,\quad \forall v_*\in V^*, \quad\forall
v\in V.
\end{equation}  
Moreover, it follows from (\ref{eq:2.24}) and (\ref{eq:2.26}) that
\begin{equation}
\label{eq:2.28}
\langle v_*,{\cal N}v_*\rangle_V = \|v_*\|^2_*\quad\,\forall\,v_*\in \mbox{dom}\,{\cal N},
\end{equation} 
and we have
\begin{equation}
\label{eq:2.29}
\langle u_*,{\cal N}v_*\rangle_V = \langle v_*,{\cal N}u_*\rangle_V =\xinto
(\nabla {\cal N} v_*)\cdot (\nabla {\cal N} u_*)\,\dx \quad\,\forall u_*,v_*\in \mbox{dom}\,{\cal N},
\end{equation}
whence also 
\begin{equation}
\label{eq:2.30}
2\,\langle \partial_t v_*(t), {\cal N}v_*(t)\rangle_V \,=\, \frac {\rm d} \dt \xinto
|\nabla{\cal N}v_*(t)|^2\,\dx \,=\,\frac {\rm d} \dt \,\|v_*(t)\|_*^2 \quad\,\mbox{for {\pier all} $t\in (0,T)$},
\end{equation}
for any $v_*\in H^1(0,T;V^*)$ satisfying $v_*^\oma(t)=0$ for a.\,e. $t\in (0,T)$.

\vspace{5mm}
The next step is to prove a priori estimates uniformly in $\alpha\in(0,1]$ for the solution $\ygesalph$ of (\ref{eq:2.13})--(\ref{eq:2.16}). 
We have the following result.

\vspace{5mm}
{\bf Proposition 2.3:}\quad\,{\em Suppose that} {\bf (A2)}--{\bf (A5)} 
{\em {\pier are} satisfied. Then there is some constant $K_1^*>0$, which only depends on $\oma$, $T$,
$y_0$, $y_{0_\Gamma}$, $f_2$, $g_2$, and $R$, such
that we have: whenever $(y^\alpha,y_\Gamma^\alpha)={\cal S}_\alpha(\uga)$ for some 
$\uga\in {\cal U}$ and some
$\alpha\in (0,1]$, then it holds}
\begin{equation}
\label{eq:2.31}
\|y^\alpha\|_{H^1(0,T;H)\cap L^\infty(0,T;V)\cap L^2(0,T;H^2(\oma))}
\,+\,\|y_\Gamma^\alpha\|_{H^1(0,T;H_\Gamma)\cap L^\infty(0,T;V_\Gamma)
\cap L^2(0,T;H^2(\Gamma))}\,\le\,K_1^*\,.
\end{equation} 

\vspace{5mm}
{\sc Proof:}\,\quad Suppose that $\uga\in {\cal U}$ and $\alpha\in (0,1]$ are arbitrarily chosen, and let $(y^\alpha,y_\Gamma^\alpha)={\cal S}_\alpha(\uga)$. The result will be established in a series of a priori
estimates. To this end, we will in the following denote by $C_i$, $i\in\nz$, positive constants which may depend on the quantities mentioned in the statement, but not on 
$\alpha\in (0,1]$. We remark that the subsequent estimates follow the same pattern as the
a priori estimates in the proof of \cite[Thm.~2.3]{CGS1}, but since not all of these 
estimates are standard, we detail them here for the reader's convenience. 

\vspace{5mm}
{\bf First a priori estimate:} \quad\,First, note that ({\pier cf.~}(\ref{eq:2.12})) $y^\alpha(t)^\oma=m_0$ for all $t\in [0,T]$,
so that $(y^\alpha(t)-m_0)\in {\rm dom}\,{\cal N}$. We thus may choose in (\ref{eq:2.17}) 
$\,v={\cal N}(y^\alpha(t)-m_0)$, and in (\ref{eq:2.18}) $\,v=-(y^\alpha(t)-m_0)$. Adding the resulting equalities,
then {\pier inserting} two additional terms on both sides for convenience, and
integrating over $[0,t]$, where $t\in [0,T]$ is arbitrary, we arrive at the identity     
\begin{align}
\label{eq:2.32}
&\frac 1 2\left(\|y^{\alpha}(t)-m_0\|^2_* + \|y^\alpha(t)-m_0\|^2_H +
\| y^{\alpha}_\Gamma(t)-m_0\|_{H_\Gamma}^2\right) \,+\,
\txt|\nabla y^\alpha|^2\,\dx\,\ds\nonumber\\[2mm]
&\quad +\txg |\nabla y_\Gamma^\alpha|^2\,\dgm\,\ds +\txt \varphi(\alpha)(h'(y^\alpha)-h'(m_0))(y^\alpha-m_0)\,\dx\,\ds
\nonumber\\[2mm]
&\quad +
\txg \psi(\alpha)(h'(y^\alpha_\Gamma)-h'(m_0))(y_\Gamma^\alpha-m_0)\,\dgm\,\ds\nonumber\\[2mm]
&=\,\frac 1 2 \left(\|y_0-m_0\|^2_* + \|y_0-m_0\|^2_H+
\| y_{0_\Gamma}-m_0\|_{H_\Gamma}^2\right) \,-\,{\Hassan \psi(\alpha)}h'(m_0)\txg (y^\alpha_\Gamma-m_0)\,\dgm\,\ds\nonumber\\[2mm]
&\quad - \txt f_2'(y^\alpha)(y^\alpha-m_0)\,\dx\,\ds + \txg (\uga-g_2'(y_\Gamma^\alpha))(y^\alpha_\Gamma -m_0)\,\dgm\,\ds\,.
\end{align}
By the monotonicity of $h'$, all of the terms on the left-hand side of (\ref{eq:2.32}) are nonnegative, while
the first term on the right-hand side is obviously bounded. Since also, in view of {\bf (A2)} and (\ref{eq:2.6}),
\begin{equation}
\label{eq:2.33}
\max_{0\le i\le 3}\left(\bigl\|f_2^{(i)}(y^\alpha)\bigr\|_{L^\infty(Q)}
\,+\,\bigl\|g_2^{(i)}(y_\Gamma^\alpha)\bigr\|_\glisig\right)\,\le\,C_1 \quad\forall\,\alpha\in (0,1], 
\end{equation}
it follows from Young's inequality and Gronwall's lemma that
\begin{equation}
\label{eq:2.34}
\|y^\alpha\|_{L^\infty(0,T;H)\cap L^2(0,T;V)} \,+\,\|y^\alpha_\Gamma\|_{L^\infty(0,T;H_\Gamma)
\cap L^2(0,T;V_\Gamma)} \,\le\,C_2\,\quad\forall\,\alpha\in (0,1].
\end{equation}

\vspace{3mm}
{\bf Second a priori estimate:} \quad\,Recalling (\ref{eq:2.12}), we may insert $\,v={\cal N}(\partial_t y^\alpha(t))$\,
in (\ref{eq:2.17}) and $\,v=-\partial_t y^\alpha(t)$ in (\ref{eq:2.18}). Adding the resulting equations,
integrating  over $[0,t]$, and using (\ref{eq:2.24}) and (\ref{eq:2.26}), we obtain the identity
\begin{align}
\label{eq:2.35}
&\int_0^t \|\partial_t y^\alpha(s)\|^2_*\,\ds \,+ \txt |\partial_t y^\alpha|^2\,\dx\,\ds\,+
\txg |\partial_t  y^{\alpha}_\Gamma|^2\,\dgm\,\ds\nonumber\\[2mm]
&\quad+\,\frac 1 2\,(\|\nabla y^\alpha(t)\|_H^2 \,+\,\|\nabla_\Gamma y^{\alpha}_\Gamma(t)\|_{H_\Gamma}^2) 
\,+\xinto \varphi(\alpha)\,h(y^\alpha(t))\,\dx\,+\ginto \psi(\alpha)\,h(y^\alpha_\Gamma(t))\,\dgm\nonumber\\[2mm]
&=\frac 1 2 \,(\|\nabla y_0\|_H^2\,+\,\|\nabla_\Gamma y_{0_\Gamma}\|^2_{H_\Gamma})
\,+\xinto \varphi(\alpha)\,h(y_0)\,\dx\,+\ginto \psi(\alpha)\,h(y_{0_\Gamma})\,\dgm
\nonumber\\[2mm]
&\quad - \txt f_2'(y^\alpha)\,\partial_t y^\alpha\,\dx\,\ds \,+\txg (\uga-g_2'(y^\alpha_\Gamma))
\,\partial _t y^\alpha_\Gamma\,\dgm\,\ds\,.
\end{align}
Obviously, the last two terms on the left-hand side 
{\pier are bounded from below and} the four terms containing the initial data
on the right-hand side of (\ref{eq:2.35}) {\pier are bounded.} 
Thus, invoking (\ref{eq:2.33}) and Young's inequality,
we can easily conclude from (\ref{eq:2.35}) the estimate
\begin{equation}
\label{eq:2.36}
\|y^\alpha\|_{H^1(0,T;H)\cap L^\infty(0,T;V)} \,+\,\|y^\alpha_\Gamma\|_{H^1(0,T;H_\Gamma)
\cap L^\infty(0,T;V_\Gamma)}\,\le\,C_3\,\quad\forall\,\alpha\in (0,1].
\end{equation}

\vspace*{3mm}
{\bf Third a priori estimate:} \quad\,Next, we insert $\,v=w^\alpha(t)-(w^\alpha(t))^\oma\,$ in (\ref{eq:2.17}) and 
apply Young's inequality, (\ref{eq:2.27}), and Poincar\'{e}'s inequality (\ref{eq:1.19}) to find the estimate
\begin{align}
\label{eq:2.37}
&\xinto \bigl|\nabla w^\alpha(t)\bigr|^2\,\dx \,=\,\xinto \bigl|\nabla \bigl(w^\alpha(t)-(w^\alpha(t))^\oma\bigr)\bigr|^2\,\dx
\,\le\,\left|\left\langle \partial_t y^\alpha(t),w^\alpha(t)-(w^\alpha(t))^\oma\right\rangle_V\right|\nonumber\\[1mm] 
&\le\,C_*\,\|\partial_t y^\alpha(t)\|_*\,\left\|w^\alpha(t)-(w^\alpha(t))^\oma\right\|_V
\,\le\,\frac 1 2 \xinto \bigl|\nabla w^\alpha(t)\bigr|^2\,\dx \,+ \,C_4\,\left\|\partial_t y^\alpha(t)\right\|_*^2\,.
%\nonumber\\[-3mm]
%&{}
\end{align}   
Now recall that the embedding $H\subset V^*$ is continuous. Hence, we can infer from estimate (\ref{eq:2.36}) that
\begin{equation}
\label{eq:2.38}
\left\|\nabla w^\alpha\right\|_{L^2(0,T;H)}\,\le\,C_5\,\quad\forall\,\alpha\in (0,1].
\end{equation}

\vspace{2mm}
Next, we aim to establish a bound for the mean value of $\,w^\alpha\,$ in $L^2(0,T)$.
To this end, we insert $v\equiv 1$ in (\ref{eq:2.18}). It follows:
\begin{align}
\label{eq:2.39}
&\xinto w^\alpha(t)\,\dx = \xinto \partial_t y^\alpha(t)\,\dx 
 +\ginto\partial_t  y^\alpha_\Gamma(t)\,\dgm + {\pier \xinto f_2'(y^\alpha(t))\,\dx }\nonumber\\[1mm] 
&\quad + {\pier \ginto (g_2'(y^\alpha_\Gamma(t))-u_\Gamma(t))\,\dgm} 
+\xinto\varphi(\alpha)h'(y^\alpha(t))\,\dx+\ginto \psi(\alpha)\,h'(y_\Gamma
^\alpha(t))\,\dgm \,.
\end{align} 
By virtue of (\ref{eq:2.33}) and (\ref{eq:2.36}), the first four integrals on the right-hand side of (\ref{eq:2.39})
define functions that are bounded in $L^2(0,T)$, uniformly in $\alpha\in (0,1]$. In order to handle the two remaining terms on the right-hand side,
we insert $v={\cal N}(y^\alpha(t)-m_0)$ in (\ref{eq:2.17}) and $v=-(y^\alpha(t)-m_0)$ in (\ref{eq:2.18})
and add the resulting equations to obtain 
\begin{eqnarray}
\label{eq:2.40}
&&\xinto\left|\nabla y^\alpha(t)\right|^2\,\dx+\ginto\left|\nabla_\Gamma y_\Gamma^\alpha(t)\right|^2\,\dgm
+\xinto \varphi(\alpha)h'(y^\alpha(t))(y^\alpha(t)-m_0)\,\dx\nonumber\\[1mm]
&&\quad +\ginto\psi(\alpha)h'(y^\alpha_\Gamma(t))(y_\Gamma^\alpha(t)-m_0)\,\dgm\,=\,G^\alpha(t), 
\end{eqnarray}
where
\begin{eqnarray}
\label{eq:2.41}
&&G^\alpha(t)\,:=\,-\xinto \partial_t y^\alpha(t)\,{\cal N}(y^\alpha(t)-m_0)\,\dx-\xinto(\partial_t y^\alpha(t)
+f_2'(y^\alpha(t)))(y^\alpha(t)-m_0)\,\dx\nonumber\\[1mm]
&&\hspace*{20mm} -\ginto \left(\partial_t y_\Gamma^\alpha(t)+g_2'(y_\Gamma^\alpha(t))-\uga(t)\right)(y^\alpha(t)-m_0)\,\dgm. 
\end{eqnarray}
Now, we may employ {\pier (\ref{eq:2.26})--(\ref{eq:2.27}) and (\ref{eq:2.33})--(\ref{eq:2.34})}  to see that
\begin{equation}
\label{eq:2.42}
|G^\alpha(t)|\,\le\,C_6\left( 1+\,\|\partial_t y^\alpha(t)\|_*\,\|y^\alpha(t)-m_0\|_*
\,+\,\|\partial_t y^\alpha(t)\|_H\,+\,\|\partial_t y_\Gamma^\alpha(t)\|_{H_\Gamma}\right),\quad
\end{equation}
for a.\,e. $t\in (0,T)$, and it follows from (\ref{eq:2.36}) that $G^\alpha$ is bounded in $L^2(0,T)$,
uniformly in $\alpha\in (0,1]$.  

\vspace{2mm}
At this point, we claim that there are $\widehat\delta>0$ and 
$\widehat C>0$ such that, for all $r\in (-1,1)$, 
\begin{equation}
\label{eq:2.43}
h'(r)(r-m_0)\,\ge\,\widehat \delta \,|h'(r)|-\widehat C\,.
\end{equation}
Indeed, since $\,-1<m_0<1$, we may employ exactly the same argument as that used in \cite[p. 908]{GiMiSchi}
to prove a corresponding estimate. From (\ref{eq:2.43}) it immediately follows that  
there is some $C_7>0$ such that for all $\alpha\in (0,1]$ {\pier we have}
\begin{align}
\label{eq:2.44}
&\varphi(\alpha)\,h'(r)(r-m_0)\ge \widehat \delta\,|\varphi(\alpha)\,h'(r)|-C_7 {\pier \quad \hbox{and}}
\nonumber\\[1mm]
& {\pier \psi(\alpha)\,h'(r)(r-m_0)\ge\widehat \delta\,|\psi(\alpha)\,h'(r)|-C_7 \quad 
\mbox{for all }\,r\in (-1,1)\,.}
\end{align}  
Consequently, we {\pier deduce that}
\begin{eqnarray}
\label{eq:2.45}
&&\xinto \varphi(\alpha)h'(y^\alpha(t))(y^\alpha(t)-m_0)\,\dx + \ginto\psi(\alpha)h'(y_\Gamma^\alpha(t))
(y_\Gamma^\alpha(t)-m_0)\,\dgm\nonumber\\[1mm]
&&\ge\,\widehat\delta \xinto |\varphi(\alpha) h'(y^\alpha(t))|\,\dx\, +\,\widehat\delta\ginto
|\psi(\alpha) h'(y_\Gamma^\alpha(t))|\,\dgm\,-\,C_8\,, 
\end{eqnarray}
and we can infer from {\pier (\ref{eq:2.39})} that
\begin{equation}
\label{eq:2.46}
\left\|(w^\alpha)^\oma\right\|_{L^2(0,T)}\,\le\,C_9 \quad\,\forall\,\alpha\in (0,1],
\end{equation}
whence, recalling (\ref{eq:2.38}) and Poincar\'e's inequality,
\begin{equation}
\label{eq:2.47}
\left\|w^\alpha\right\|_{L^2(0,T;V)}\,\le\,C_{10}\quad\,\forall\alpha\in (0,1]\,.
\end{equation}

\vspace{5mm}
{\bf Fourth a priori estimate:} \quad\,{\pier Next,} observe that in view of (\ref{eq:2.19}), {\Hassan (\ref{eq:2.20})}
and (\ref{eq:2.22}) we have $(v,v_\Gamma)\in {\cal V}$ for $v=\varphi(\alpha)h'(y^\alpha)$.
Hence, we may insert $v=\varphi(\alpha)h'(y^\alpha)$ in (\ref{eq:2.18}) to
obtain
\begin{eqnarray}
\label{eq:2.48}
&&\txt \varphi(\alpha)\,h''(y^\alpha)\left|\nabla y^\alpha\right|^2\,\dx\,\ds +
\txg \varphi(\alpha)\,h''(y^\alpha_\Gamma)\left|\nabla_\Gamma y^\alpha_\Gamma
\right|^2\,\dgm\,\ds \nonumber\\[1mm] 
&&\quad +\txt |\varphi(\alpha)\,h'(y^\alpha)|^2\,\dx\,\ds +
\txg \varphi(\alpha)\,\psi(\alpha)\,|h'(y^\alpha_\Gamma)|^2\,\dgm\,\ds \nonumber\\[1mm] 
&&=\txt \varphi(\alpha)\,h'(y^\alpha)(w^\alpha-f_2'(y^\alpha)-\partial_t y^\alpha)\,
\dx\,\ds\nonumber\\[1mm]
&&\quad + \txg \varphi(\alpha)\,h'(y^\alpha_\Gamma)(\uga-g_2'(y^\alpha_\Gamma)
-\partial_t y^\alpha_\Gamma)\,\dgm\,\ds\,.
\end{eqnarray}
{\Hassan Now notice that $\,h''>0\,$ in $(-1,1)$, which implies that the two integrals in which $\,h''\,$ occurs in the integrands,
are both nonnegative. Moreover, (\ref{eq:1.12}) implies that
\begin{equation*}
\txg \varphi(\alpha)\psi(\alpha)\,|h'(y_\Gamma^\alpha)|^2\,\dgm\,\ds \,\ge\,
\frac 1 {C_{\varphi\psi}}\txg (\varphi(\alpha))^2\,|h'(y_\Gamma^\alpha)|^2\,\dgm\,\ds\,.
\end{equation*}
Therefore the boundary integral 
\begin{align*}
\txg \varphi(\alpha)\,h'(y^\alpha_\Gamma)(\uga-g_2'(y^\alpha_\Gamma)
-\partial_t y^\alpha_\Gamma)\,\dgm\,\ds
\end{align*}
can be handled using Young's inequality. Now applying (\ref{eq:2.33}), (\ref{eq:2.36}), (\ref{eq:2.47}) and Young's inequality, we find that}
\begin{equation}
\label{eq:2.49}
\left\|\varphi(\alpha)h'(y^\alpha)\right\|_{L^2(0,T;H)}\,\le\,C_{11} \quad\,
\forall\,\alpha\in (0,1]\,.
\end{equation}  

\vspace{5mm}
{\bf Fifth a priori estimate:} \quad\,Now observe that the variational equality
(\ref{eq:2.18}) implies that $y^\alpha$ solves (\ref{eq:2.14}) at least in the sense of
distributions. Since all other terms have been proved to be bounded in $L^2(0,T;H)$,
we must have
\begin{equation}
\label{eq:2.50}
\|\Delta y^\alpha\|_{L^2(0,T;H)}\,\le\,C_{12} \quad\,\forall\,\alpha\in (0,1]\,.
\end{equation} 
Next, we use \cite[Thm.~3.2, p.~1.79]{Brezzi} to conclude that 
$$
\int_0^T\|y^\alpha(t)\|^2_{H^{3/2}(\oma)}\,\dt \le \,C_{13}
\int_0^T (\|\Delta y^\alpha(t)\|_H^2 + {\pier \|y^\alpha_\Gamma (t)\|_{V_\Gamma}^2} )\,\dt\,,
$$
whence it follows that
\begin{equation}
\label{eq:2.51}
\|y^\alpha\|_{L^2(0,T;H^{3/2}(\oma))}\,\le\,C_{14}\,\quad\forall\,\alpha\in
(0,1]\,.
\end{equation}
Hence, by the trace theorem \cite[Thm.~2.27, p.~1.64]{Brezzi}, we have
\begin{equation}
\label{eq:2.52}
\|\partial_\nf y^\alpha\|_{L^2(0,T;H_\Gamma)}\,\le\,C_{15}\,\quad\forall\,\alpha\in
(0,1]\,.
\end{equation}

From the above estimates it follows that all the terms occurring in the integration
by parts formula for the Laplace operator are functions, and we deduce that the 
variational equation (\ref{eq:2.18}) also implies that the second identity in
(\ref{eq:2.15}) holds at least in a generalized sense, in principle. Therefore, the
preceding estimates yield that, {\Hassan {\pier by letting} $G_\Gamma^\alpha:=\uga \,-\, \partial_\nf y^\alpha {}-{}
\partial_t y_{\Gamma}^\alpha \,-\,  g_{2}'(y_\Gamma^\alpha) $}, we can write
\begin{equation}
\label{eq:2.53}
-\Delta_\Gamma y^\alpha_\Gamma+\psi(\alpha)h'(y^\alpha_\Gamma)={\Hassan G_\Gamma^\alpha}
\,\,\mbox{on $\Sigma$, where }\,\|{\Hassan G_\Gamma^\alpha}\|_{L^2(\Sigma)}
\,\le\,C_{16} \,\,\,\forall\,\alpha \in (0,1]\,.
\end{equation}
Testing the above equation by $\,\psi(\alpha)h'(y^\alpha_\Gamma)$, we obtain
\begin{eqnarray}
\label{eq:2.54}
&&\txg \psi(\alpha)h''(y^\alpha_\Gamma)\,\left|\nabla_\Gamma y_\Gamma^\alpha
\right|^2\,\dgm\,\ds+ \txg |\psi(\alpha)h'(y^\alpha_\Gamma)|^2\,\dgm\,\ds\nonumber\\[1mm]
&&=\txg \psi(\alpha)h'(y_\Gamma^\alpha)\,{\Hassan G_\Gamma^\alpha} \,\dgm\,\ds\,, 
\end{eqnarray}
and {\Hassan a simple application of Young's inequality shows that} 
\begin{equation}
\label{eq:2.55}
\|\psi(\alpha)h'(y^\alpha_\Gamma)\|_{L^2(0,T;H_\Gamma)}\,\le\,C_{17}\,\quad\forall\,
\alpha\in (0,1],
\end{equation} 
whence also 
\begin{equation}
\label{eq:2.56}
\|\Delta_\Gamma y^\alpha_\Gamma\|_{L^2(0,T;H_\Gamma)}\,\le\,C_{18}\,\quad\forall\,
\alpha\in (0,1]\,.
\end{equation}
The boundary version of the elliptic regularity theory then yields
\begin{equation}
\label{eq:2.57}
\|y^\alpha_\Gamma\|_{L^2(0,T;H^2(\Gamma))}\,\le\,C_{19}\,\quad\forall\,
\alpha\in (0,1]\,,
\end{equation} 
and {\pier consequently} it follows from standard elliptic estimates that
\begin{equation}
\label{eq:2.58}
\|y^\alpha\|_{L^2(0,T;H^2(\oma))}\,\le\,C_{20}\,\quad\forall\,
\alpha\in (0,1]\,.
\end{equation}
With this, the assertion is completely proved.\qed

\section{Existence and approximation of optimal controls}\label{existence}
\setcounter{equation}{0}

{Our first aim in this section is to prove the following existence result:

\vspace{5mm}
{\bf Theorem 3.1:}\quad\,{\em Suppose that the assumptions} {\bf (A1)}--{\Hassan {\bf (A5)}} 
{\em are satisfied. Then the optimal control problem} $({\cal P}_0)$ {\em admits a solution.}
 
\vspace{5mm}
Before proving Theorem~3.1, we introduce the solution space
\begin{eqnarray}
\label{eq:3.1}
&&{\cal Y}:= \left\{(y,\yga)\in {\cal V}: y\in H^1(0,T;H)\cap
L^\infty(0,T;V)\cap L^2(0,T;H^2(\oma)),\right.\qquad\nonumber\\[1mm]
&&\qquad\quad\left. y_\Gamma=y_{{\pier |_\Gamma}}, \quad  y_\Gamma\in H^1(0,T;H_\Gamma)\cap
L^\infty(0,T;V_\Gamma)\cap L^2(0,T;H^2(\Gamma))\right\}\,,
\end{eqnarray}
and a family of auxiliary optimal control problems $({\cal P}_\alpha)$, which is
parametrized by $\,\alpha \in(0,1]\,$. In what follows, we will always assume that $h$ is given by {(\ref{eq:1.10})} and that $\,\varphi\,$ and $\,\psi\,$ are functions that are positive and continuous on $(0,1]$ and satisfy the conditions (\ref{eq:1.12}). 
For $\,\alpha \in(0,1]$, let us denote by ${\mathcal S}_\alpha$ the operator mapping 
$\uga\in\uad$ into the unique solution $(y^\alpha,y_\Gamma^\alpha) \in {\cal Y}$  to the  
variational problem {\pier (\ref{eq:2.16})}--(\ref{eq:2.18}). We define: 
\begin{eqnarray*}
({\mathcal{P}}_{\alpha})\quad&&\mbox{Minimize }\,\, {\cal J}((y,\yga),\uga)
\quad\mbox{over }\,{\mathcal Y}\times \uad \quad\mbox{subject to the condition that}
\nonumber\\  &&\mbox{{\pier (\ref{eq:2.16})}--(\ref{eq:2.18}) are satisfied.}
\end{eqnarray*}

The following result is a consequence of \cite[{\pier Thm.~2.2}]{CGS2}.

\vspace{3mm}
{\bf Lemma~3.2:}\,\quad {\em Suppose that the assumptions} {\bf (A1)}--{\Hassan {\bf (A5)}}
{\em and} {\rm (\ref{eq:1.10}), (\ref{eq:1.12})}
{\em are fulfilled, and let $\alpha\in(0,1]$ be given. Then the optimal control 
problem $({\mathcal{P}}_{\alpha})$ admits a solution.}

\vspace{7mm}
{\sc Proof of Theorem~3.1:} \,\quad
%\vspace{2mm} 
Let $\,\{\alpha_n\}\subset (0,1]\,$ be any sequence such that $\alpha_n\searrow 0$ as $n\to\infty$. By virtue of 
Lemma~3.2, {\pier for any $n\in\nz$ we may} pick an optimal pair for the
optimal control problem $({\cal P}_{\alpha_n})$,
$${((y^{\alpha_n},y^{\alpha_n}_\Gamma),u^{\alpha_n}_\Gamma)}\in{\cal Y}\times\uad$$ 
where $(y^{\alpha_n},y_\Gamma^{\alpha_n},w^{\alpha_n})$ is the unique solution to {\pier (\ref{eq:2.16})}--(\ref{eq:2.18}),
written for $\alpha=\alpha_n$,
which satisfies (\ref{eq:2.19})--(\ref{eq:2.22}). In particular, 
$\,(y^{\alpha_n},y^{\alpha_n}_\Gamma)\,=\,{\cal S}_{\alpha_n}(u^{\alpha_n}_\Gamma)$ for all $n\in\nz$. 
Moreover, Proposition {\Hassan 2.3}  implies that (\ref{eq:2.31}) holds for any $\alpha_n$, $n\in\nz$. From this and
from (\ref{eq:2.47}) we may without loss of generality assume that there are $\,\uga\in\uad$, $w$, and $\,(y,y_\Gamma)\,$  such that
\begin{eqnarray}
\label{eq:3.2}
&&u^{\alpha_n}_\Gamma\to\uga\quad\mbox{weakly-star in }\,{\cal X}\,,\\[1mm]
\label{eq:3.3}
&&w^{\alpha_n}\to w\quad\mbox{weakly in } \,L^2(0,T;V)\,,\\[1mm]
\label{eq:3.4}
&&y^{\alpha_n}\to y\quad\mbox{weakly-star in } \,H^1(0,T;H)\cap L^\infty(0,T;V)\cap L^2(0,T;H^2(\Omega))\,,\\[1mm]
\label{eq:3.5}
&&{y^{\alpha_n}_\Gamma}\to y_\Gamma\quad\mbox{weakly-star in }\,H^1(0,T;H_\Gamma)\cap L^\infty(0,T;V_\Gamma)\cap L^2(0,T;H^2(\Gamma))\,.\quad 
\end{eqnarray} 
By the continuity of the embedding 
$H^1(0,T;H) \cap L^2(0,T;H^2(\oma){\Hassan )}\subset C^0([0,T];V)$,
we have in fact $y\in C^0([0,T];V)$, and, by the same token, $y_\Gamma\in C^0([0,T];V_\Gamma)$. 
Owing to the Aubin-Lions lemma {(see \cite[Sect.~8, Cor.~4]{Simon})}, we also have
\begin{eqnarray}
\label{eq:3.6}
&&{y^{\alpha_n}}\to  y\quad\mbox{strongly in }\,C^0([0,T];H)\cap L^2(0,T;V)\,, \\[1mm]
\label{eq:3.7}
&&{y^{\alpha_n}_\Gamma}\to  y_\Gamma\quad\mbox{strongly in }\,
C^0([0,T];H_\Gamma)\cap L^2(0,T;V_{\Gamma}). 
\end{eqnarray}
In particular, it holds $\,y(\cdot,0)=y_0\,$, as well as 
$\, y_\Gamma(\cdot,0)=y_{0_\Gamma}$.  
In addition, the Lipschitz continuity of $f'_2$ and $g'_2$ on $[-1,1]$ yields that
\begin{eqnarray}
\label{eq:3.8}
&&{f'_2({y^{\alpha_n}})\to f'_2(y)}\quad\mbox{strongly in }\,C^0([0,T];H),\\[1mm]
\label{eq:3.9}
&&{g'_2({y^{\alpha_n}_\Gamma})\to g'_2(y_\Gamma)}\quad\mbox{strongly in }\,C^0([0,T];H_\Gamma)\,. 
\end{eqnarray}
Moreover, (\ref{eq:2.49}) and (\ref{eq:2.55}) show that without loss of generality we 
may also assume that
\begin{eqnarray}
\label{eq:3.10}
&&\varphi(\alpha_n)\,h'({y^{\alpha_n}})\to \xi \quad\mbox{weakly in }\,L^2(0,T;H),
\\[1mm]
\label{eq:3.11}
&&\psi(\alpha_n)\,h'({y^{\alpha_n}_\Gamma})\to \xi_\Gamma\quad\mbox{weakly in }\,L^2(0,T;H_\Gamma), 
\end{eqnarray}
for some weak limits $\xi $ and $\xi_\Gamma$. 

\vspace{2mm}
Combining the above convergences, we may pass to the limit
as $n\to\infty$ in (\ref{eq:2.17}) and (\ref{eq:2.18}) (written for $\alpha_n$) to find that the
quintuple $(y,\yga, w, \xi,\xiga)$ is a solution to (\ref{eq:2.10})--(\ref{eq:2.11}), and obviously 
the properties (\ref{eq:2.4})--(\ref{eq:2.6}) and (\ref{eq:2.9}) are satisfied. In order to show that 
the quintuple $(y,\yga,w,\xi,\xiga)$ is a solution to problem (\ref{eq:1.3})--(\ref{eq:1.7}) in the sense
of Definition 2.1, it remains to
 show that $\,\xi\in\partial I_{[-1,1]}(y)\,$ a.\,e. in $Q$
and $\,\xiga\in \partial I_{[-1,1]}(y_\Gamma)\, $ a.\,e. in $\Sigma$. 
Once this will be shown, we can conclude 
that $(y, y_\Gamma) = {\cal S}_0
(\uga)$, i.\,e., that the pair $\,(( y,  y_\Gamma,w,\xi,\xiga),\uga)\,$ is admissible for $({\cal P}_0)$. 

\vspace{2mm}
Now, recalling (\ref{eq:1.10}) and owing to the convexity of $\,h$, we have, for every $n\in\nz$,
\begin{align}
\label{eq:3.12}
&\txinto\varphi(\alpha_n)\,h(y^{\alpha_n})\,\dx\,\dt\,+\,\txinto\varphi(\alpha_n)\,h'( y^{\alpha_n})\,{{} (z-y^{\alpha_n})}\,\dx\,\dt\,\leq\, \txinto\varphi(\alpha_n)\,h(z)\,\dx\,\dt\nonumber\\[2mm]
&\quad \mbox{for all }\,z\in {\cal K}=\{v\in {L^2(Q)}:|v|\leq 1\text{ a.e. in }Q\}\,. 
\end{align}

Thanks to (\ref{eq:1.12})}, the integral on the right-hand side and the first integral on the left-hand side of 
(\ref{eq:3.12}) tend to zero as $n\to\infty${\pier , since $h$ is a bounded function}. Hence, {invoking (\ref{eq:3.6}) and (\ref{eq:3.10}), {\pier the} passage to the limit as $n\to\infty$} yields
\begin{equation}
\label{eq:3.13}
\txinto\xi\,( y-z)\,\dx\,\dt\,\geq 0\quad\forall z\in {\mathcal{K}} . 
\end{equation}
Inequality (\ref{eq:3.13}) entails that $\xi$ is an element of the subdifferential of the extension $\mathcal{I} $ of $ I_{[-1,1]}$ to $L^2(Q)$, which means that $\xi \in \partial \mathcal{I} (y)$ or, equivalently (cf.~\cite[Ex.~2.3.3., p.~25]{Brezis}),  
$\xi\in\partial I_{[-1,1]}(y)$ {a.\,e. in $Q$}. Similarly we prove that 
$\xi_\Gamma\in\partial I_{[-1,1]}(y_\Gamma)$ {a.\,e. in $\Sigma$}.

\vspace{2mm}
It remains to show that $((y,y_\Gamma,w,\xi,\xiga),\uga)$ is in fact optimal for $({\cal P}_0)$.
To this end, let $\vga\in \uad$ be arbitrary. In view of the convergence properties (\ref{eq:3.2}) and
{(\ref{eq:3.4})--(\ref{eq:3.7})},
and using the weak sequential lower semicontinuity properties of the cost functional, we have
\begin{align}
\label{eq:3.14}
&{\cal J}((y,y_\Gamma),\uga)\,=\,{\cal J}({\cal S}_{0}(\uga),\uga)\,\le\,
\liminf_{n\to\infty}\,{\cal J}({\mathcal S}_{\alpha_n}(u^{\alpha_n}_\Gamma),u^{\alpha_n}_\Gamma)
 \nonumber\\[1mm]
& {\pier \,\leq\,\liminf_{n\to\infty}\,{\cal J}({\cal S}_{\alpha_n}(v_{\Gamma}),v_{\Gamma}) =}\lim_{n\to\infty} {\cal J}({\cal S}_{\alpha_n}(v_{\Gamma}),v_{\Gamma})\,=\,
{\cal J}({\mathcal S}_{0}(v_{\Gamma}),v_{\Gamma}),
\end{align}  
where for the last equality the continuity of the cost functional with respect to the first variable
was used. With this, the assertion is completely proved.\qed

\vspace{5mm}
{\bf Corollary 3.3:}\quad\,{\em Let the general assumptions} {\bf (A1)}--{\bf (A5)} {\em and}
{\rm (\ref{eq:1.10}), (\ref{eq:1.12})}  {\em be satisfied, and
let sequences {$\, \{\alpha_n\}\subset (0,1]\,$ and $\,\{u^{\alpha_n}_\Gamma\}\subset {\cal U}$}  be given such that,
as $n\to\infty$, $\,\alpha_n\searrow 0\,$ and $\,u^{\alpha_n}_\Gamma\to u_\Gamma\,$ weakly-star in $\,{\cal X}$.
Then we have}
\begin{eqnarray}
\label{eq:3.15}
&&{\cal S}_{\alpha_n}(u^{\alpha_n}_\Gamma)\to\,{\cal S}_0 (u_\Gamma)\quad\mbox{\em weakly-star in }\,{\cal Y}\,,
\\[2mm]
\label{eq:3.16}
&&\lim_{n\to\infty} {\cal J}({\cal S}_{\alpha_n}(v_\Gamma),v_\Gamma)\,=\,{\cal J}({\cal S}_0(\vga),
\vga) \quad\forall \,\vga\in {\cal U}\,.
\end{eqnarray}

\vspace{3mm}
{\sc Proof:}\quad\,By the same arguments as in the first part
of the proof of Theorem~3.1, we can conclude that 
(\ref{eq:3.15}) holds at least for some subsequence. But the limit is given by the first two
components of a solution quintuple in the sense of Definition 2.1 to the state system 
(\ref{eq:1.3})--(\ref{eq:1.7}), which, according to Proposition 2.2, are uniquely determined. 
Hence, the limit is the same for all convergent subsequences and 
(\ref{eq:3.15}) is true for the entire sequence. Now, let $\vga\in {\cal U}$ be arbitrary.
Then (see (\ref{eq:3.6})--(\ref{eq:3.7}))
$\,{\cal S}_{\alpha_n}(\vga)\,$ converges strongly to $\,{\cal S}_0(\vga)\,$ in 
$\,(C^0([0,T];H)\cap L^2(0,T;V))\times (C^0([0,T];H_\Gamma)\cap L^2(0,T;V_{\Gamma}))$,
so that (\ref{eq:3.16}) follows from the continuity properties of the cost functional with respect to
its first argument.\qed

\vspace{7mm}
Theorem~3.1 does not yield any information on whether every solution to the optimal control problem $({\mathcal{P}}_{0})$ can be approximated by a sequence of solutions to the problems $({\mathcal{P}}_{\alpha})$. 
As already announced in the Introduction, we are not able to prove such a general `global' result. Instead, we 
can only give a `local' answer for every individual optimizer of $({\mathcal{P}}_{0})$. For this purpose,
we employ a trick due to Barbu \cite{Barbu}. To this end, let $\bar u_\Gamma\in\uad$
be an arbitrary optimal control for $({\mathcal{P}}_{0})$, and let $(\bar y,\bar y_\Gamma,\bar w,\bar\xi,\bar\xiga)$
be an associated solution quintuple to the state system (\ref{eq:1.3})--(\ref{eq:1.7}) in the sense of 
Definition 2.1. In particular, $\,(\bar y,\bar y_\Gamma)={\cal S}_0 (\bar u_\Gamma)$. We associate with this 
optimal control the {\em adapted cost functional}
\begin{equation}
\label{eq:3.17}
\widetilde{\cal J}((y,\yga),\uga):={\cal J}((y,\yga),\uga)\,+\,\frac{1}{2}\,\|u_{\Gamma}-\bar{u}_\Gamma\|^2_{L^2(\Sigma)}
\end{equation}
and a corresponding {\em adapted optimal control problem}
\begin{eqnarray*}
(\widetilde{\mathcal{P}}_{\alpha})\qquad\mbox{Minimize }\,\, \widetilde {\cal J}((y,\yga),\uga)
\quad\mbox{over }\,{\mathcal Y}\times \uad \quad\mbox{subject to the condition}\hspace*{17mm}\nonumber\\  \mbox{that (\ref{eq:2.13})--(\ref{eq:2.16}) be satisfied.}\hspace*{81.5mm}
\end{eqnarray*}

\vspace{3mm}
With a standard direct argument that needs no repetition here, we can show the following 
result.

\vspace{5mm}
{\bf Lemma 3.4:}\quad\,{\em Suppose that the assumptions} {\bf (A1)}--{\Hassan {\bf (A5)}}
{\em and} {\rm (\ref{eq:1.10}), (\ref{eq:1.12})} {\em are
satisfied, and let $\alpha\in (0,1]$. Then the optimal control problem} 
$(\widetilde{\cal P}_\alpha)$ {\em admits a solution.}
  
\vspace{5mm}
We are now in the position to give a partial answer to the question raised above. We have the following result.

\vspace{5mm}
{\bf Theorem~3.5:}\,\quad{\em Let the general assumptions} {\bf (A1)}--{\Hassan {\bf (A5)}} 
{\em and}  {\rm (\ref{eq:1.10}), (\ref{eq:1.12})} {\em be fulfilled, and suppose that 
$\bar u_\Gamma\in \uad$ is an arbitrary optimal control of} $({\mathcal{P}}_{0})$ {\em with associated state
quintuple $(\bar y,\bar y_\Gamma, \bar w, \bar\xi,\bar\xiga)$. 
Then for every sequence $\,\{\alpha_n\}\subset (0,1]$ such that
$\,\alpha_n\searrow 0\,$ as $\,n\to\infty$ and for any $n\in\nz$ there exists some optimal control
 $\,\bar u_\Gamma^{\alpha_n}\in \uad\,$ of the adapted problem $(\widetilde{\mathcal{P}}_{\alpha_n})$
 with associated state triple $(\bar y^{\alpha_n},\bar y_\Gamma^{\alpha_n}, \bar w^{\alpha_n})$ such that, as $n\to\infty$,}
\begin{eqnarray}
\label{eq:3.18}
\hskip-.8cm &&\bar u_\Gamma^{\alpha_n}\to \bar u_\Gamma\quad\mbox{strongly in }\,H_\Gamma,\\[1mm]
\label{eq:3.19}
\hskip-.8cm &&\bar y^{\alpha_n} \to \bar y \quad\mbox{weakly-star in }\,H^1(0,T;H)\cap L^\infty(0,T;V)\cap L^2(0,T;H^2(\oma)),\\[1mm]
\label{eq:3.20}
\hskip-.8cm && \bar y_\Gamma^{\alpha_n}\to \bar{y}_\Gamma\quad\mbox{weakly-star in } \,
H^1(0,T;H_\Gamma)\cap L^\infty(0,T;V_\Gamma)\cap L^2(0,T;H^2(\Gamma))\,, \qquad\\[1mm]
\label{eq:3.21}
\hskip-.8cm &&\widetilde{{\cal J}}((\bar y^{\alpha_n},\bar y_\Gamma^{\alpha_n}),\bar u_\Gamma^{\alpha_n})\to  {\cal J}((\bar{y},\bar{y}_\Gamma),\bar u_\Gamma)\,.
\end{eqnarray}

\vspace{3mm}
{\sc Proof:} \quad\,Let $\alpha_n \searrow 0$ as $n\to\infty$. For any $ n\in\nz$, we pick an optimal control $\bar u_\Gamma^{\alpha_n} \in \uad\,$ for the adapted problem $(\widetilde{\cal P}_\alpha)$ and denote by 
$(\bar y^{\alpha_n},\bar y_\Gamma^{\alpha_n}, \bar w^{\alpha_n})$ the associated solution triple of 
problem (\ref{eq:2.13})--(\ref{eq:2.16}); in particular, we have $(\bar y^{\alpha_n},\bar y_\Gamma^{\alpha_n})
={\cal S}_{\alpha_n}(\bar u_\Gamma^{\alpha_n})$, and (\ref{eq:2.19})--(\ref{eq:2.22}) are satisfied.
By the boundedness of $\uad$, we have for some subsequence of $\{\alpha_n\}$, which is again indexed by $n$,
that  it holds  
\begin{equation}
\label{eq:3.22}
\bar u_\Gamma^{\alpha_n}\to \uga\quad\mbox{weakly-star in }\,{\cal X}
\quad\mbox{as }\,n\to\infty,
\end{equation}
with some $\uga\in\uad$. Owing to Corollary 3.3, we have
\begin{equation}
\label{eq:3.23}
(\bar y ^{\alpha_n},\bar y^{\alpha_n}_\Gamma )={\cal S}_{\alpha_n}(\bar u_\Gamma^{\alpha_n})\to {\cal S}_0(\uga)=:(y,\yga) \quad 
\mbox{weakly-star in }\,{\cal Y}.
\end{equation}
In particular, $y,\yga$ are the first two components of a quintuple $(y,\yga,w,\xi,\xiga)$ solving 
the state system associated with $\uga$, which implies that $((y,\yga,w,\xi,\xiga),\uga)$
is admissible for $({\cal P}_0)$.

\vspace{2mm}
We now aim to prove that $\uga=\bar u_\Gamma$. Once this will be shown, the uniqueness result of Proposition 2.2 yields 
that also $\,(y,\yga)=(\bar y,\bar y_\Gamma)$, which shows that (\ref{eq:3.19}) and (\ref{eq:3.20}) 
hold at least for the subsequence; but since the limit is the same for any subsequence, we
have (\ref{eq:3.19}),
(\ref{eq:3.20}) for the entire sequence $\{\alpha_n\}$. By the same token, also
(\ref{eq:3.22}) will hold for the entire
sequence.  

\vspace{2mm}
Indeed, we have, owing to the weak sequential lower semicontinuity of $\widetilde {\cal J}$, 
and in view of the optimality property
of  $\,((\bar{y},\bar{y}_\Gamma),\buuga)$ for problem $({\cal P}_0)$,
\begin{eqnarray}
\label{eq:3.24}
&&\liminf_{n\to\infty}\, \widetilde{\cal J}((\bar y^{\alpha_n},\bar y_\Gamma^{\alpha_n}),\bar u_\Gamma^{\alpha_n})
\ge \,{\cal J}((y,\yga),\uga)\,+\,\frac{1}{2}\,
\|u_{\Gamma}-\bar{u}_\Gamma\|^2_{L^2(\Sigma)}\nonumber\\[1mm]
&&\geq \, {\cal J}((\bar{y},\bar{y}_\Gamma),\bar u_\Gamma)\,+\,\frac{1}{2}\,\|u_{\Gamma}-\bar{u}_\Gamma\|^2_{L^2(\Sigma)}\,.
\end{eqnarray}
On the other hand, the optimality property of  $\,((\bar{y}^{\alpha_n},\bar{y}_\Gamma^{\alpha_n}),
\bar u_\Gamma^{\alpha_n})\,$ for problem $(\widetilde {\cal P}_{\alpha_n})$ yields that
for any $n\in\nz$ we have
\begin{equation}
\label{eq:3.25}
\widetilde {\cal J}((\bar{y}^{\alpha_n},\bar{y}_\Gamma^{\alpha_n}),
\bar u_\Gamma^{\alpha_n})\, =\,
\widetilde {\cal J}({\cal S}_{\alpha_n}(\bar{u}_\Gamma^{\alpha_n}),
\bar u_\Gamma^{\alpha_n})\,\le\,\widetilde {\cal J}({\cal S}_{\alpha_n}(\bar u_\Gamma),\bar u_\Gamma)\,,
\end{equation}
whence, taking the limes superior as $n\to\infty$ on both sides and invoking (\ref{eq:3.16}) in
Corollary~3.3,
\begin{eqnarray}
\label{eq:3.26}
&&\limsup_{n\to\infty}\,\widetilde {\cal J}((\bar{y}^{\alpha_n},\bar{y}_\Gamma^{\alpha_n}),
\bar u_\Gamma^{\alpha_n})\,\le\,\widetilde {\cal J}({\cal S}_0(\bar u_\Gamma),\bar u_\Gamma) \,=\,
\widetilde {\cal J}((\bar y,\bar y_\Gamma),\bar u_\Gamma)\qquad\quad\nonumber\\[1mm]
&&=\,{\cal J}((\bar y,\bar y_\Gamma),\buuga)\,.
\end{eqnarray}
Combining (\ref{eq:3.24}) with (\ref{eq:3.26}), we have thus shown that 
$\,\frac{1}{2}\,\|u_{\Gamma}-\bar{u}_\Gamma\|^2_{L^2(\Sigma)}=0$\,,
so that $\,\uga=\buuga\,$  and thus also $\,(y,\yga)=(\bar y,\bar y_\Gamma)$.
Moreover, (\ref{eq:3.24}) and (\ref{eq:3.26}) also imply that
\begin{eqnarray}
\label{eq:3.27}
&&{\cal J}((\bar{y},\bar{y}_\Gamma),\buuga) \, =\,\widetilde{\cal J}((\bar{y},\bar{y}_\Gamma),\buuga)
\,=\,\liminf_{n\to\infty}\, \widetilde{\cal J}((\bar y^{\alpha_n},\bar y_\Gamma^{\alpha_n}),
\bar u_\Gamma^{\alpha_n})\nonumber\\[1mm]
&&\,=\,\limsup_{n\to\infty}\, \widetilde{\cal J}((\bar y^{\alpha_n},\bar y_\Gamma^{\alpha_n}),  \bar u_\Gamma^{\alpha_n})\,
=\,\lim_{n\to\infty}\, \widetilde{\cal J}((\bar y^{\alpha_n},\bar y_\Gamma^{\alpha_n}),
\bar u_\Gamma^{\alpha_n})\,,
\end{eqnarray}                                     
which proves {(\ref{eq:3.21})} and, at the same time, also (\ref{eq:3.18}). The assertion is thus
completely {\pier checked}.\qed

\section{The optimality system}\label{optimality}
\setcounter{equation}{0}
In this section our aim is to establish first-order necessary optimality conditions for the optimal control problem $({\mathcal{P}}_{0})$.  This will be achieved by passage to the limit as $\alpha\searrow 0$ in the {\pier (recently 
in \cite{CGS2})} derived first-order necessary optimality conditions for the adapted optimal control problems $(\widetilde{\mathcal{P}}_{\alpha})$. It will turn out that in the limit certain generalized first-order necessary conditions of optimality result. To fix things once and for all,
we will throughout the entire section assume that $h$ is given by (\ref{eq:1.10}) and that (\ref{eq:1.12}) and
the general assumptions {\bf (A1)}--{\bf (A5)} are
satisfied; we also assume that a fixed optimal control $\bar u_\Gamma\in \uad$ for $({\cal P}_0)$, along
with a solution quintuple $(\bar y,\bar y_\Gamma, \bar w, \bar\xi,\bar\xiga)$ of the associated state system
(\ref{eq:1.3})--(\ref{eq:1.7}), is given. 
In additon, we make the following compatibility assumption:

\vspace{5mm}
{\bf (A6)} \quad\,It holds $\beta_3=\beta_4=0$.

\vspace{5mm} 
We remark that in \cite[Remark 5.6]{CGS2} it has been pointed out that this assumption it dispensable
at the expense of less regularity of the adjoint state variables; in order to keep the technicalities
at a reasonable level, we here confine {\pier ourselves} to the case $\beta_3=\beta_4=0$. 

\subsection{The optimality conditions for $(\widetilde{\mathcal{P}}_{\alpha})$}
 We begin our analysis by formulating the adjoint state system for the adapted
control problem $(\widetilde{\mathcal{P}}_{\alpha})$.
To this end, let us assume that $\bar u_\Gamma^\alpha\in\uad$ is an arbitrary optimal control for 
$(\widetilde{\mathcal{P}}_{\alpha})$ and that $(\bar y^\alpha,\bar y_\Gamma^\alpha, \bar w^\alpha)$ 
is the solution triple to the associated state system (\ref{eq:2.13})--(\ref{eq:2.16}). In particular, 
$(\bar y^\alpha, \bar y_\Gamma^\alpha)={\cal S}_\alpha(\bar u_\Gamma^\alpha)$, and the solution has the 
regularity properties (\ref{eq:2.19})--(\ref{eq:2.22}). It then follows (see \cite[Eqs. (5.7)--(5.9)]{CGS2})
that the corresponding adjoint state variables $\,q^\alpha,q^\alpha_\Gamma,p^\alpha\,$ solve the following backward-in-time variational problem:
\begin{align}
\label{eq:4.1}
&\xinto q^\alpha(t)\,v\,\dx\,=\,\xinto \nabla p^\alpha(t)\cdot \nabla v\,\dx \quad\,{\pier \hbox{ for all }\,v\in V \,\,
\mbox{and }}\,t\in (0,T)\,,\\[2mm]
\label{eq:4.2}
&-\xinto \partial_t\left(q^\alpha(t)+p^\alpha(t)\right)\,v\,\dx + \xinto\nabla q^\alpha(t)\cdot\nabla v\,\dx
+\ginto\nabla_\Gamma q_\Gamma^\alpha(t)\cdot \nabla_\Gamma \vga\,\dgm\nonumber\\[1mm]
&\quad -\ginto \partial_t q_\Gamma^\alpha \,\vga\,\dgm + \xinto \bigl(\varphi(\alpha)h''(\bar y^\alpha(t))
 +f_2''(\bar y^\alpha(t))\bigr)q^\alpha(t)\,v\,\dx\nonumber\\[1mm]
&\quad + \ginto \bigl(\psi(\alpha)h''(\bar y^\alpha_\Gamma(t))
 +g_2''(\bar y^\alpha_\Gamma(t))\bigr)q^\alpha_\Gamma(t)\,\vga\,\dgm\nonumber\\[1mm]
&=\xinto \beta_1(\bar y^\alpha(t)-z_Q(t))\,v\,\dx + \ginto \beta_2(\bar y_\Gamma^\alpha(t)-z_\Sigma(t))\,\vga
\,\dgm\nonumber\\[1mm]
&\quad\,{\pier \hbox{for every }}\,(v,\vga)\in {\cal V} \quad\mbox{and a.\,a. }\,t\in (0,T)\,,
\\[2mm]
\label{eq:4.3}
&\xinto (q^\alpha(T)+p^\alpha(T))\,v\,\dx + \ginto q^\alpha(T)\,\vga \,\dgm =0 \quad\,{\pier \hbox{for every }}\,
(v,\vga)\in {\cal V}\,.
\end{align}

In \cite[Thm.~2.4]{CGS2} it has been shown that the system (\ref{eq:4.1})--(\ref{eq:4.2})
has for every $\alpha \in (0,1]$ a unique solution triple $(q^\alpha,q^\alpha_\Gamma,
p^\alpha)$ such that
\begin{equation}
\label{eq:4.4}
(q^\alpha,q^\alpha_\Gamma)\in {\cal Y},\,\quad 
p^\alpha \in H^1(0,T;H^2(\oma))\cap  L^2(0,T;H^4(\oma)),
\end{equation}
and we may regard $(q^\alpha,q^\alpha_\Gamma,p^\alpha)$ as a solution to the 
linear PDE system
\begin{align}
\label{eq:4.5}
&-\Delta p^\alpha=q^\alpha \quad\mbox{in }\,Q,\,\quad\,\partial_\nf p^\alpha=0 
\quad\mbox{on }\,\Sigma,\\[2mm]
\label{eq:4.6}
&-\partial_t (q^\alpha+p^\alpha)-\Delta q^\alpha+\bigl(\varphi(\alpha)h''(\bar
y^\alpha)+f_2''(\bar y^\alpha)\bigr)q^\alpha=\beta_1(\bar y^\alpha-z_Q)
\quad\mbox{in }\,Q,\quad\\[2mm]
\label{eq:4.7}
&{\pier {}- \partial_t q_\Gamma^\alpha}+\partial_\nf q^\alpha-\Delta_\Gamma q_\Gamma^\alpha
+\bigl(\psi(\alpha)h''(\bar y_\Gamma^\alpha)+g_2''(\bar y_\Gamma^\alpha)\bigr)
q_\Gamma^\alpha=\beta_2(\bar y_\Gamma^\alpha-z_\Sigma)
\nonumber\\[1mm]
&{\pier \quad\mbox{and}\quad} q^\alpha_{{\pier |_\Gamma}}=q^\alpha_\Gamma \,\quad\mbox{on }\,\Sigma,\quad\\[2mm]
\label{eq:4.8}
&q^\alpha(T)+p^\alpha(T)=0 \quad\mbox{in }\,\oma,\,\quad\,q^\alpha_\Gamma(T)=0
\quad\mbox{on } \,\Gamma.
\end{align}
Moreover, {\pier as we are now dealing with 
$(\widetilde{\mathcal{P}}_{\alpha})$ instead of 
$({\mathcal{P}}_{\alpha})$, the variational inequality given by 
\cite[Thm.~2.5]{CGS2} has to be modified as follows:}
\begin{equation}
\label{eq:4.9}
\int_0^T\!\!\ginto \bigl(q^\alpha_\Gamma + \beta_5\,\bar u_\Gamma^\alpha
+(\bar u_\Gamma^\alpha-\bar u_\Gamma)\bigr)(\vga-\bar u_\Gamma^\alpha)\,\dgm\,\dt \,\ge\,0 \quad
\forall\,\vga \in \uad\,.
\end{equation}
 
In order to pave the road for the limit process as $\alpha\searrow 0$ in the 
optimality conditions for $(\widetilde{\mathcal{P}}_{\alpha})$, we employ an idea
that was developed in \cite{CGS2}. Namely, it is possible to show that the system
(\ref{eq:4.1})--(\ref{eq:4.3}) is equivalent to a decoupled problem that can
be solved by first finding $q^\alpha$ and then reconstructing $p^\alpha$. We briefly 
motivate this approach. First, standard embedding results yield that
$q^\alpha\in C^0([0,T];V)$, and it
immediately follows from inserting $v\equiv 1$ in (\ref{eq:4.1}) that 
$(q^\alpha({\Hassan t}))^\oma=0$ for all $t\in [0,T]$. Hence $q^\alpha(t)\in {\rm dom}\,
{\cal N}$, and, with the mean value function $(p^\alpha)^\oma \in C^0([0,T])$,
the function $(p^\alpha-(p^\alpha)^\oma)(t)$ satisfies for every $t\in [0,T]$
the identity (\ref{eq:2.24}) with $v_*=q^\alpha(t)$. In other words, we have
\begin{equation}
\label{eq:4.10}
(p^\alpha-(p^\alpha)^\oma)(t)={\cal N}(q^\alpha(t))\quad\forall\,t\in [0,T].
\end{equation}
On the other hand, $(p^\alpha(t))^\oma$ is for any fixed $t\in [0,T]$ a constant
function and thus orthogonal in $H$ to the subspace of functions having zero
mean value. Consequently, $p^\alpha$ is completely eliminated from (\ref{eq:4.2})
if we confine ourselves to the use of test functions having zero mean value.
Similar remarks apply for the final condition on $q^\alpha+p^\alpha$ appearing
in (\ref{eq:4.3}). In this way, we may try to first construct $(q^\alpha,q^\alpha
_\Gamma)$ and then recover $p^\alpha$ from (\ref{eq:4.10}), where the calculation
of $(p^\alpha(t))^\oma$ is an easy task,
since simple integration of (\ref{eq:4.6}) over $\oma\times [t,T]$, using (\ref{eq:4.8})
and the fact that $q^\alpha(t)$ has zero mean value, immediately yields that
$$
(p^\alpha(t))^\oma\,=\,\frac 1 {|\oma|}
\int_t^T\!\!\xinto \left(-\Delta q^\alpha
+ \bigl(\varphi(\alpha) h''(\bar y^\alpha)+f_2''(\bar y^\alpha)\bigr)\,
q^\alpha - \beta_1(\bar y^\alpha-z_Q)\right)\,\dx\,\ds\,.
$$

\vspace{2mm}
We now make this approach precise. Since our test functions
will have zero mean value, we introduce the linear spaces
\begin{equation}
\label{eq:4.11}
{\cal H}_\oma:=\left\{(v,\vga)\in{\cal H}: v^\oma=0\right\}, \quad\,
{\cal V}_\oma:= {\cal H}_\oma\cap {\cal V},
\end{equation}
and we define on ${\cal H}_\oma$ and ${\cal V}_\oma$ the inner products
\begin{eqnarray}
\label{eq:4.12} 
&&((u,\uga),(v,\vga))_{{\cal H}_\oma}\,:=\,((u,\uga),(v,\vga))_{{\cal H}}\,=\,\xinto u\,v\,\dx+
\ginto \uga \,\vga\,\dgm, \qquad\\[2mm]
\label{eq:4.13}
&&((u,\uga),(v,\vga))_{{\cal V}_\Omega}\,:=\,\xinto\nabla u\cdot\nabla v\,\dx +\ginto
\nabla_\Gamma \uga\cdot \nabla_\Gamma \vga\,\dgm,
\end{eqnarray}
where $(u,\uga),(v,\vga)$ are generic elements of ${\cal H}_\oma$ (resp., ${\cal V}_\oma$). 
Note that it follows from Poincar\'{e}'s inequality (\ref{eq:1.19}) that (\ref{eq:4.13}) actually
defines an inner product in ${\cal V}_\oma$ whose associated norm is equivalent to the
standard one. 

\vspace{2mm}
Next, we infer from \cite[Lemma 5.1 and Cor. 5.3]{CGS2} that 
\begin{equation}
\label{eq:4.14}
V_\Gamma=\left\{\vga: (v,\vga)\in {\cal V}_\oma\right\}, \quad\mbox{and ${\cal V}_\oma$
is dense in ${\cal H}_\oma$}. 
\end{equation}
Therefore, we can construct the Hilbert triple ${\cal V}_\oma\subset {\cal H}_\oma\subset {\cal V}_\oma^*$ with
dense and compact embeddings,
that is, we identify ${\cal H}_\oma$ with a subspace of ${\cal V}_\oma^*$ in such a way that
\begin{equation}
\label{eq:4.15}
\langle (u,\uga),(v,\vga)\rangle_{{\cal V}_\oma}\,=\,((u,\uga),(v,\vga))_{{\cal H}_\oma}
\quad\,\forall\,(u,\uga)\in {\cal H}_\oma, \quad\forall\,(v,\vga)\in {\cal V}_\oma\,.
\end{equation}

\vspace{2mm}
Observe that{\pier , because of the zero mean value condition,} the first components $v$
of the elements $(v,\vga)\in {\cal V}_\oma$ cannot span the whole space 
$C_0^\infty(\oma)$; consequently, variational {\pier equalities} with test functions in
${\cal V}_\oma$ cannot immediately be interpreted as equations in the sense of 
distributions. We obviously have the following result:

\vspace{5mm}
{\bf Lemma 4.1:} \quad\,{\em Let the general assumptions} {\bf (A1)}--{\bf (A6)}
{\em and} {\rm (\ref{eq:1.10}), (\ref{eq:1.12})} {\em be satisfied. Then the
pair $(q,q_\Gamma)=(q^\alpha,q^\alpha_\Gamma)$ is a solution to the 
variational system}
\begin{align}
\label{eq:4.16}
&-\xinto\partial_t\bigl({\cal N}(q(t))+q(t)\bigr)v\,\dx+\xinto\nabla q(t)\cdot
\nabla v\,\dx+ \ginto \nabla_\Gamma q_\Gamma(t)\cdot\nabla_\Gamma \vga\,\dgm
\nonumber\\[1mm]
&\quad -\ginto \partial_t q_\Gamma\,\vga\,\dgm +\xinto \bigl(\varphi(\alpha)
h''(\bar y^\alpha(t))+f_2''(\bar y^\alpha(t))\bigr)\,q(t)\,v\,\dx\nonumber\\[1mm]
&\quad +\ginto \bigl(\psi(\alpha)h''(\bar y_\Gamma^\alpha(t))+g_2''(\bar y_\Gamma
^\alpha(t))\bigr)\,q_\Gamma(t)\,\vga\,\dgm\nonumber\\[1mm]
&=\xinto \beta_1(\bar y^\alpha(t)-z_Q(t))\,v\,\dx + \ginto\beta_2
(\bar y^\alpha_\Gamma(t)-z_\Sigma(t))\,\vga\,\dgm\nonumber\\[1mm]
&\quad{\pier \hbox{for every }} \,(v,\vga)\in {\cal V}_\oma \quad\mbox{{\em and for a.\,a.} }\,
t\in (0,T),\\[2mm]
\label{eq:4.17}
&\xinto {\Hassan ({\cal N}(q)+q)(T)}\,v\,\dx +\ginto q_\Gamma(T)\,\vga\,\dgm=0\,\quad
{\pier \hbox{for every }} \,(v,\vga)\in {\cal V}_\oma\,.
\end{align}

\vspace{5mm}
Notice that we may insert $(v,\vga)=(q^\alpha(T),q^\alpha_\Gamma(T))\in {\cal V}_\oma$
in the end point condition (\ref{eq:4.17}), which, in view of (\ref{eq:2.28}), yields
that
$$
{\pier \|q^\alpha(T)\|_*^2\,+\, \|q^\alpha(T)\|_H^2}\,+\,\|q^\alpha_\Gamma(T)\|_{H_\Gamma}^2
\,=\,0\,;
$$
we thus may replace (\ref{eq:4.17}) by the simpler condition
\begin{equation}
\label{eq:4.18}
q^\alpha(T)=0 \quad\mbox{a.\,e. in }\,\oma, \,\quad\,{\cal N}(q^\alpha(T))=0 \quad\mbox{a.\,e. in }\,\oma, \,\quad\,q^\alpha_\Gamma(T)=0\quad
\mbox{a.\,e. on }\,\Gamma ,
\end{equation}
where the second equation {\pier simply follows 
from the fact that $q^\alpha(T)$ belongs to the domain of 
the operator~${\cal N}$}.
% from the fact that $z={\cal N}(q^\alpha(T))$ has zero 
% mean value and solves the Neumann problem \,\,$-\Delta 
% z=q^\alpha(T)=0\,$ in $\,\oma$, \,$\partial_\nf z=0\,$ on $\,\Gamma$. 

\vspace{5mm}
{\bf Remark 4.2:} \quad\,In \cite[{\pier Theorems~2.5 and~5.4}]{CGS2} it has been shown that there is only one solution to problem (\ref{eq:4.16})--(\ref{eq:4.17}) (namely,
$(q^\alpha,q^\alpha_\Gamma)$) that has zero mean value and belongs to ${\cal Y}$. 

\vspace{5mm}
We now prove an a priori estimate which will be fundamental for the derivation of the optimality conditions for $(\mathcal{P}_{0})$. To this end, we introduce some further function spaces. At first, we put
\begin{equation}
\label{eq:4.19}
{\pier {\mathcal W}}\,:= {\pcgg \bigl( H^1(0,T; V^*) \times H^1(0,T; V_\Gamma^*) \bigr)}
\cap L^2(0,T;{\cal V}_\oma).
\end{equation}
Then we define 
\begin{equation}
\label{eq:4.20}
{\pier {\mathcal W}_0} \,:\,=\{(\eta,\eta_{\Gamma})\in {\pier {\mathcal W}} :(\eta {(0)},\eta_{\Gamma}{(0))}=(0,0)\}\,.
\end{equation}
Observe that both these spaces are Banach spaces when equipped with the natural norm of $\,
{\pier {\cal W}} $. Moreover, ${\pier {\mathcal W}} $ is continuously embedded in
{\pcgg $C^0([0,T];H) \times C^0([0,T]; H_\Gamma)$}, 
so that the initial condition encoded in (\ref{eq:4.20}) is 
meaningful.
{\pcgg Furthermore, since ${\mathcal W}_0$ is a closed subspace of 
\begin{equation}
\bigl( H^1(0,T; V^*) \cap L^2(0,T;V)\bigr)  
\times \bigl( H^1(0,T; V_\Gamma^*) \cap L^2(0,T;V_\Gamma)\bigr) ,\label{pier1}
\end{equation}
we deduce that the elements $F \in {\mathcal W}_0^* $ are exactly the ones given by 
\begin{align}
\label{rappr}
\left\langle \!\left\langle F,(\eta,\eta_\Gamma)\right\rangle\! \right\rangle \,=\,
\left\langle z, \eta  \right\rangle  +
\left\langle z_\Gamma,\eta_\Gamma \right\rangle_\Gamma \quad 
\hbox{ for all $(\eta,\eta_\Gamma)\in {\pier {\mathcal W}_{0}}$,}
\end{align}
where $z$ and $z_\Gamma$ vary in the dual spaces of $ H^1(0,T; V^*) \cap L^2(0,T;V)$ 
and $ H^1(0,T; V_\Gamma^*) \cap L^2(0,T;V_\Gamma)$, respectively. 
Of course, the duality symbols in (\ref{rappr}) refer to ${\mathcal W}_0 $ and the two spaces 
above and their corresponding duals; moreover, the representation of $F$ through 
(\ref{rappr}) is not unique. Notice that a particular $z$ might be any function in $L^2(0,T; V^*)$ by means of the natural embedding 
\begin{equation}
 L^2(0,T; V^*) \subset 
\bigl( H^1(0,T; V^*) \cap L^2(0,T;V)\bigr)^* \nonumber
\end{equation}
(due to the density of $H^1(0,T; V^*) \cap L^2(0,T;V) $ in $ L^2(0,T; V)$), i.e.,  
$$
\left\langle z, v  \right\rangle  = \int_0^T \left\langle z(t) , v(t)  \right\rangle_V dt  
\quad \hbox{ for all } v\in H^1(0,T; V^*) \cap L^2(0,T;V).
$$
Analogously, we can take $z_\Gamma \in L^2(0,T;V_\Gamma^*)$.
Finally, the above representation formula allows us to give a meaning to a sentence like 
$$
(z^\alpha , z_\Gamma^\alpha ) \to  (z , z_\Gamma )
\quad \hbox{ weakly in } {\mathcal W}_0^* .  
$$
}%
Next, we put
\begin{equation}
\label{eq:4.22}
{\mathcal Z}\,:=\,{\pcgg L^\infty\left(0,T;{\cal H}_\oma\right)}\cap L^2\left(0,T;{\cal V}\right),
\end{equation}
which is a Banach space when equipped with its natural norm.

\vspace{5mm}
{\bf Proposition 4.3:} \quad\,{\em Let the general assumptions} {\bf (A1)}--{\bf (A6)}
{\em and} {\rm (\ref{eq:1.10}), (\ref{eq:1.12})} {\em be satisfied and let}
\begin{equation}
\label{eq:4.23}
(\lambda^\alpha,\lambda_\Gamma^\alpha)\,:=\,\left(\varphi(\alpha)\,h''(\bar{y}^{\alpha})\,q^{\alpha}\,, \, \psi(\alpha)\,h''(\bar{y}_\Gamma^{\alpha})\,q_\Gamma^{\alpha}\right)
\quad\forall\,\alpha\in (0,1].
\end{equation}
{\em Then there exists
a constant $K_2^*>0$, which only depends on the data of the system and on $R$, such
that for all $\alpha\in (0,1]$ it holds}
\begin{align}
\label{eq:4.24}
&{\pcgg \left\|(q^\alpha,q_\Gamma^\alpha)\right\|_{\cal Z}\,+\,
\left\|(\lambda^{\alpha},\lambda^{\alpha}_\Gamma)\right\|_{\pier {\cal W}_0^*}\,+\,
\left\|{\cal N}(q^\alpha)\right\|_{L^\infty(0,T;H^2(\oma))\cap L^2(0,T;H^3(\oma))}}
\nonumber\\[2mm]
% &\,+\,\left\|\left({\cal N}(q^\alpha)\right)_{{\pier |_\Gamma}}\right\|_{L^\infty(0,T;
% H^{3/2}(\Gamma))\cap L^2(0,T;H^{5/2}(\Gamma))}
% \nonumber\\[2mm]
&{\pcgg \, +\,\left\|
(\partial_{t}\left({\cal N}(q^\alpha)+q^{\alpha}\right),
\partial_{t}{\pcgg q_\Gamma^{\alpha}})\right\|_{\pier {\cal W}_0^*} 
\,\le\,K_2^*\,.}
\end{align}

\vspace{3mm}
{\sc Proof:} \quad\,In the following, $C_i$, $i\in\nz$, denote positive constants which are independent
of $\alpha\in (0,1]$. To show the boundedness of the adjoint variables, we insert $(v,\vga)=(q^\alpha(t),
q^\alpha_\Gamma(t))\in {\cal V}_\oma$ in (\ref{eq:4.16}), written for $(q,q_\Gamma)=(q^\alpha,q^\alpha_\Gamma)$,
and integrate over $[s,T]$ where $s\in [0,T]$. {\pier First, note that 
\begin{align}
\label{eq:4.25}
&-\int_s^T\!\!\xinto\partial_t\left({\cal N}(q^\alpha)+q^\alpha\right)\,q^\alpha\,\dx\,\dt
=\xinto\left({\cal N}(q^\alpha(s))q^\alpha(s)+\frac 1 2\,\left|q^\alpha(s)\right|^2\right)\dx\nonumber\\[2mm]
&\quad+\int_s^T\!\!\xinto {\cal N}(q^\alpha)\,\partial_t q^\alpha\,\dx\,\dt \,=\,\frac 1 2 \left
(\left\|q^\alpha(s)\right\|^2_H\,+\,\left\|q^\alpha(s)\right\|_*^2\right)
\end{align} 
since $\,\partial_t{\cal N}(q^\alpha)\in L^2(0,T;H^2(\oma))$ by (\ref{eq:4.4}) and (\ref{eq:4.10}), 
and the integration by parts with respect to time can be done in view of  
(\ref{eq:2.28}), (\ref{eq:2.30}), and (\ref{eq:4.18}).}
We thus obtain the equation
\begin{align}
\label{eq:4.26}
&\frac{1}{2}\left(\|q^\alpha(s)\|^2_H\,+\,\|q^\alpha(s)\|^2_*\,+\,\|q^\alpha_\Gamma(s)\|^2_{H_\Gamma}\right)
+\int_s^T\!\!\xinto\left|\nabla q^\alpha\right|^2\dx\,\dt   \nonumber\\[2mm]
&\quad {\pier {}+ \int_s^T\!\!\ginto
\left|\nabla_\Gamma q^\alpha_\Gamma\right|^2\dgm\,\dt }+\int_s^T\!\!\xinto \lambda^\alpha\,q^\alpha\,\dx\,\dt +
\int_s^T\!\!\xinto \lambda^\alpha_\Gamma\,q^\alpha_\Gamma\,\dgm\,\dt\nonumber\\[2mm]
&{\pier {}= - \int_s^T\!\!\xinto f_2''(\bar{y}^{\alpha})\left|q^{\alpha}\right|^2\, \dx\, \dt\,{\pcgg {}-{}}\int_s^T\!\!\ginto g_2''(\bar{y}_\Gamma^{\alpha})\left|q_\Gamma^{\alpha}\right|^2\, \dgm\, \dt}\nonumber\\[2mm]
&\quad +\int_s^T\!\!\xinto\beta_1\,\left(\bar{y}^{\alpha}-z_Q\right)\,q^{\alpha}\,\dx\,\dt\,+\,
\int_s^T\!\!\ginto\beta_2\left(\bar{y}_\Gamma^{\alpha}-z_\Sigma\right)\,q_\Gamma^
{\alpha}\,\dgm\,\dt\,.
\end{align}

\vspace{2mm}
By (\ref{eq:4.23}) and the positivity of $h''$, the {\pier last} 
two integrals in the second line of the left-hand 
side of (\ref{eq:4.26}) are nonnegative, while, {\pier owing to (\ref{eq:2.33}) and {\bf (A1)},
the right-hand side} of (\ref{eq:4.26}) can obviously be bounded by an
expression of  the form
$$
C_1\,\Bigl(1\,+\int_s^T\!\!\xinto\left|q^{\alpha}\right|^2\, \dx\, \dt\,
+\,\int_s^T\!\!\ginto \left|q_\Gamma^{\alpha}\right|^2\, \dgm\, \dt\Bigr)\,.
$$
Hence, invoking Gronwall's inequality, we find the estimate
\begin{equation}
\label{eq:4.27}
\|(q^{\alpha},q^\alpha_\Gamma)\|_{L^\infty\left(0,T;{\cal H}\right)\cap L^2\left(0,T;{\cal V}\right)}
\,\leq C_2\,\quad\forall\,\alpha\in (0,1]\,.
\end{equation}
Moreover, {\pcgg using (\ref{eq:2.25}) we find that
\begin{align}
\label{eq:4.28}
&\|{\cal N}(q^\alpha)\|_{L^\infty(0,T;H^2(\oma))\cap L^2(0,T;H^3(\oma))}\,\le\,C_3
\,\quad\forall\,\alpha\in (0,1].
%\\[2mm]
% \label{eq:4.29}
% &\left\|({\cal N}(q^\alpha))_{{\pier |_\Gamma}}\right\|_{L^\infty(0,T;H^{3/2}(\Gamma))\cap
% L^2(0,T;H^{5/2}(\Gamma))}\,\le\,C_4\,\quad\forall\,\alpha \in (0,1]\,.
\end{align}
}%

\vspace{3mm}
Next, we derive the bound for the time derivatives. To this end, let 
$\,(\eta,\eta_\Gamma)\in{\pier {\mathcal W}_0} \,$ be arbitrary. 
Using (\ref{eq:4.18}),  the initial condition for $(\eta,\eta_\Gamma)$, and the 
estimates (\ref{eq:4.27})--(\ref{eq:4.28}), we obtain from integration by parts  
% {\pier and (\ref{eq:4.21})} 
that
\begin{align}
\label{eq:4.30}
&\left\langle\! \left\langle-(\partial_t ({\cal N}(q^\alpha)+q^\alpha),
\partial_t {\pcgg q_\Gamma^{\alpha}})\,,\,(\eta,\eta_{\Gamma})\right\rangle\!\right\rangle\nonumber\\[2mm] 
&=-\txinto \partial_t\left({\cal N}(q^\alpha)+q^\alpha\right)\eta\,\dx\,\dt\,-\,\tgamma\partial_t {\pcgg q_\Gamma^{\alpha}}
\,\eta_\Gamma \,\dgm\,\dt\nonumber\\[2mm]
&=\tinto\!\!\left\langle \partial_t \eta(t),{\cal N}(q^\alpha(t))+
q^\alpha(t)\right\rangle_V\dt\,+\,\tinto\!\!\left\langle\partial_t\eta_\Gamma(t),
{\pcgg q_\Gamma^{\alpha}}(t)\right\rangle_{V_\Gamma}\dt \nonumber\\[2mm]
&\le \tinto\|\partial_t\eta(t)\|_{V^*}\,\|{\cal N}(q^\alpha(t))+q^\alpha(t)\|_V
\,\dt + \tinto \|\partial_t\eta_\Gamma(t)\|_{V_\Gamma^*}\,\left\|
{\pcgg q_\Gamma^{\alpha}}(t)\right\|_{V_\Gamma}\dt\nonumber\\[2mm]
&\le {\pcgg C_4} \,\|(\eta, \eta_\Gamma)\|_{{\pcgg {\cal W}_0}}\,, \quad\mbox{for all }\,
\alpha\in (0,1]\,.
\end{align}
We thus have shown that
\begin{equation}
\label{eq:4.31}
\left\|
(\partial_{t}({\cal N}(q^\alpha)+q^\alpha),\partial_t{\pcgg q_\Gamma^{\alpha}}
)\right\|_{\pier {\cal W}_0^*}
\,\le {\pier \,{\pcgg C_4}}\,\quad\forall \,\alpha\in (0,1]\,.
\end{equation}
Finally, {\pier by recalling (\ref{eq:4.23}) and the estimates (\ref{eq:4.27})--(\ref{eq:4.28}), (\ref{eq:4.31}),   
a comparison in (\ref{eq:4.16}) yields that
\begin{equation}
\label{eq:4.32}
\|(\lambda^{\alpha},\lambda_\Gamma^{\alpha})\|_{\pier {\cal W}_0^*}\,\le\,{\pcgg C_5}\,\quad\forall\,\alpha\in (0,1]
\end{equation}
as well,} and the assertion is proved.
\qed

\subsection{The optimality conditions for $({\mathcal{P}}_{0})$}.  

We now establish first-order necessary optimality conditions for $({\mathcal{P}}_0)$
by performing a limit as $\alpha\searrow 0$ in the approximating problems. To this end, recall that a fixed optimal control $\bar u_\Gamma\in \uad$ for $({\cal P}_0)$, along
with a solution quintuple $(\bar y,\bar y_\Gamma, \bar w, \bar\xi,\bar\xiga)$ of the associated state system (\ref{eq:1.3})--(\ref{eq:1.7}) is given.

\vspace{2mm}
We draw some consequences from the previously established results. First recall that
by Theorem~3.5 for any sequence $\{\alpha_n\}\subset (0,1]$ with $\alpha_n\searrow 0$ as
$n\to\infty$, and for any $n\in\nz$ we can find an optimal control $\bar u_\Gamma^{\alpha_n}\in\uad$ for
$(\widetilde{\mathcal{P}}_{\alpha_n})$ and an associated state triple $(\bar y^{\alpha_n},\bar y_\Gamma^{\alpha_n},\bar w^{\alpha_n})$ such that the convergences (\ref{eq:3.18})--(\ref{eq:3.21}) hold. As in the proof of Theorem~3.1, we may without loss of generality assume
that 
\begin{eqnarray}
\label{eq:4.33}
&&f_2''(\bar y^{\alpha_n})\to f_2''(\bar y)\,\quad\mbox{strongly in }\,C^0([0,T];H),\\[2mm]
\label{eq:4.34}
&&g_2''(\bar y^{\alpha_n}_\Gamma)\to g_2''(\bar y_\Gamma)\,\quad\mbox{strongly in }\,C^0([0,T];H_\Gamma)\,.
\end{eqnarray}
Also, by virtue of {\Hassan Lemma 4.1} and Proposition 4.3, we may without loss of generality assume
that there exist the corresponding {\Hassan adjoint state variables $(q^{\alpha_n},
q^{\alpha_n}_\Gamma)\in {\cal Y}$ that satisfy} {\pcgg
\begin{align}
\label{eq:4.35}
&(q^{\alpha_n},q_\Gamma^{\alpha_n})\to (q,q_\Gamma)\quad\mbox{weakly-star in }\,
{\cal Z},\\[2mm]
\label{eq:4.36}
&{\cal N}(q^{\alpha_n}) \to {\cal N}(q)\quad\mbox{weakly-star in }\,
L^\infty(0,T;H^2(\oma))\cap L^2(0,T;H^3(\oma)),\\[2mm] 
% \label{eq:4.37} 
%&({\cal N}(q^{\alpha_n}))_{{\pier |_\Gamma}}\to ({\cal N}(q))_{{\pier |_\Gamma}} \quad\mbox{weakly-star in }
%\,L^\infty(0,T;H^{3/2}(\Gamma))\cap {\pier L^2}(0,T;H^{5/2}(\Gamma)), \\[2mm]
\label{eq:4.38}
&(\lambda^{\alpha_n},\lambda_\Gamma^{\alpha_n})\to (\lambda,\lambda_\Gamma)
\quad\mbox{weakly in } {\cal W}_0^*\,,
\end{align}
}%
for suitable limits $\,(q,q_\Gamma)\,$ and \,$(\lambda,\lambda_\Gamma)${\pcgg , where 
$\lambda$ and $\lambda_\Gamma$ belong to the duals of the spaces involved in (\ref{pier1}),
as explained above.}
Therefore, passing to the limit as
$n\to\infty$ in the variational inequality (\ref{eq:4.9}), written for $\alpha_n$, $n\in\nz$, {\pier and recalling  
(\ref{eq:3.18}),} we obtain that $(q,q_\Gamma)$ satisfies
\begin{equation}
\label{eq:4.39}
\tgamma \!(q_\Gamma\,+\,\beta_5\,{\bar u}_\Gamma)\,(v_\Gamma-
{\bar u}_\Gamma)\,\dgm\,\dt\,\ge\,0 \quad
\forall\,v_\Gamma\in\uad.
\end{equation}

\vspace*{2mm}
Next, we will show that in the limit as $n\to \infty$ a limiting adjoint system for $({\cal P}_0)$
is satisfied. To this end, we insert an arbitrary $(\eta,\eta_\Gamma)\in {\cal W}_0\,$  in (\ref{eq:4.16}),
written for $\alpha_n$, $n\in\nz$, and 
integrate the resulting equation over $[0,T]$. Integrating by parts with respect to $t$, and invoking (\ref{eq:4.18})
and the zero initial conditions for $(\eta,\eta_\Gamma)$,  we arrive at the identity
\begin{align}
\label{eq:4.40}
&\txinto\lambda^{\alpha_n}\,\eta\,\dx\,\dt\,+\!\tgamma \lambda_\Gamma^{\alpha_n}\,\eta_\Gamma\,\dgm\,\dt
\,+\!\int_0^T\!\langle {\partial_t} \eta (t),{\cal N}(q^{\alpha_n}(t))+q^{\alpha_n}(t)\rangle_V\,\dt\nonumber\\[2mm]
&\quad +\tinto\!\langle
\partial_t\eta_\Gamma(t),q_\Gamma^{\alpha_n}(t)\rangle_{V_\Gamma}\,\dt
 +\txinto \nabla q^{\alpha_n}\cdot \nabla\eta\,\dx\,\dt
\,+\,\tgamma \nabla_\Gamma q_\Gamma^{\alpha_n}\cdot\nabla_\Gamma \eta_\Gamma\,\dgm\,\dt\nonumber\\[2mm]
&\quad+\txinto f_2''({\bar y^{\alpha_n}})\,q^{\alpha_n}\,\eta\,\dx\,\dt \,+\,\tgamma g_2''({\bar y^{\alpha_n}_\Gamma})\,
q_\Gamma^{\alpha_n}\,\eta_\Gamma\,\dgm\,\dt\nonumber\\[2mm]
&= \beta_1\txinto ({\bar y^{\alpha_n}}-z_Q)\,\eta\,\dx\,\dt\,+\,\beta_2\tgamma ({\bar y^{\alpha_n}_\Gamma}-z_\Sigma)\,\eta_\Gamma
\,\dgm\,\dt\,.
\end{align}

Now, by virtue of the convergences (\ref{eq:3.19}), (\ref{eq:3.20}), and (\ref{eq:4.33})--(\ref{eq:4.38}), we may pass to the limit as
$n\to\infty$ in (\ref{eq:4.40}) to obtain, for all $ \,(\eta,\eta_\Gamma)\in {\cal W}_0$,		
\begin{align}
\label{eq:4.41}
&\langle\!\langle (\lambda,\lambda_\Gamma),(\eta,\eta_\Gamma)\rangle\!\rangle \,+\,
\int_0^T\!\langle {\partial_t} \eta (t)\,,\, {\cal N}(q(t))+q(t)\rangle_V\,\dt\,+\tinto\!\langle
\partial_t\eta_\Gamma(t)\,,\, q_\Gamma(t)\rangle_{V_\Gamma}\,\dt\nonumber\\[2mm]
&\quad+\txinto\nabla q\cdot\nabla\eta\,\dx\,\dt\,+\,\tgamma\nabla_\Gamma q_\Gamma\cdot
\nabla_\Gamma \eta_\Gamma\,\dgm\,\dt\nonumber\\[2mm]
&\quad+\txinto f_2''(\bar y)\, q\,\eta\,\dx\,\dt \,+\,\tgamma g_2''(\bar y_\Gamma)\,
q_\Gamma\,\eta_\Gamma\,\dgm\,\dt\nonumber\\[2mm]
&=\beta_1\txinto (\bar y-z_Q)\,\eta\,\dx\,\dt\,+\,\beta_2\tgamma (\bar y_\Gamma-z_\Sigma)\,\eta_\Gamma
\,\dgm\,\dt\,.
\end{align}

\vspace*{3mm}
Next, we show that the limit pair $\,((\lambda,\lambda_\Gamma),(q,q_\Gamma))\,$ satisfies some sort
of a complementarity slackness condition. To this end, observe that for all $n\in\nz$ we obviously
have
$$
\txinto\lambda^{\alpha_n}\,q^{\alpha_n}\,\dx\,\dt\,=\,\txinto\varphi(\alpha_n)\,h''(\bar y^{\alpha_n})
\,|q^{\alpha_n}|^2\,\dx\,\dt\,\ge\,0\,.$$
An analogous inequality holds for the corresponding boundary terms. We thus have
\begin{equation}
\label{eq:4.42}
{\liminf_{n\to\infty} \txinto\lambda^{\alpha_n}\,q^{\alpha_n}\,\dx\,\dt\, \ge 0,\quad \liminf_{n\to\infty} \tgamma\lambda_\Gamma^{\alpha_n}\,q_\Gamma^{\alpha_n}\,
\dgm\,\dt\,\ge\,0\,.}
\end{equation}
   
\vspace*{3mm}
Finally, we derive a relation which gives some indication that the limit $(\lambda,\lambda_\Gamma)$ 
should somehow be  concentrated
on the set where $\,|\bar y|=1\,$ and $\,|\bar y_\Gamma|=1$ (which, however, we cannot prove rigorously). 
To this end, we test the pair $\,(\lambda^{\alpha_n},\lambda_\Gamma^{\alpha_n})\,$  by the function
$\,\left((1-(\bar y^{\alpha_n})^2)\,\phi, (1-(\bar y_\Gamma^{\alpha_n})^2)\,\phi_\Gamma\right)\,$ {\pier that belongs to
${\cal V}_\Omega$ since} $\,(\phi,\phi_\Gamma)\,$ is any smooth test function satisfying 
\begin{equation}
\label{eq:4.43}
(\phi(0),\phi_\Gamma(0))
=(0,0), \quad \int_\oma (1-(\bar y^{\alpha_n}(t))^2)\,\phi(t)\,\dx=0 \quad\forall\,t\in [0,T].
\end{equation}
{\pier As $\,h''(r)= 2 / \left( 1-r^2 \right) $ for every $r\in (-1,1)$,} we obtain
\begin{align}
\label{eq:4.44}
&\lim_{n\to\infty}\left(\txinto\lambda^{\alpha_n}\,(1-(\bar{y}^{\alpha_n})^2)\,\phi\,\dx\,\dt\,,\tgamma\lambda_\Gamma^{\alpha_n}\,(1-(\bar{y}_\Gamma^{\alpha_n})^2)\,\phi_\Gamma
\,\dgm\,\dt\right)\nonumber\\[2mm]
&=\,\lim_{n\to\infty}\left(2\txinto \varphi(\alpha_n)\,q^{\alpha_n}\,\phi\,\dx\,\dt\,,\,
2\tgamma\psi(\alpha_n)\, q_\Gamma^{\alpha_n}\,\phi_\Gamma\,\dgm\,\dt\right)\,=\,(0,{0})\,.\quad
\end{align}

\vspace{5mm}
We now collect the results established above, {especially} in Theorem~3.5. We have the following statement.

\vspace{7mm}
{\bf Theorem~4.4:}\,\quad{\em Let the assumptions} {\bf (A1)}--{\bf (A6)} {\em be satisfied, let $h$ be given by}
(\ref{eq:1.10}), {\em and let $\varphi,\psi$ be positive and continuous functions on $(0,1]$ fulfilling} (\ref{eq:1.12}).
{\em Moreover, let $\,\bar u_\Gamma\in\uad$ be an optimal control for $({\cal P}_0)$ with
associated solution quintuple $(\bar y,\bar y_\Gamma, \bar w,\bar\xi,\bar\xiga)$ to the corresponding state system} (\ref{eq:1.3})--(\ref{eq:1.7}) {\em in the sense of Definition 2.1.
Then the following assertions hold true:}

\vspace{2mm}
(i) \,\,\,{\em For every sequence} $\{\alpha_n\}\subset (0,1]$, 
{\em with} $\,\alpha_n\searrow 0\,$ {\em as} 
$n\to\infty$, {\em and for any} $n\in\nz$, {\em there exists a solution} 
$\,\bar u_\Gamma^{\alpha_n}\in\uad\,$ {\em to the adapted control
problem} $\,(\widetilde{\mathcal{P}}_{\alpha_n})$ {\em such that, with the associated solution
triple} $(\bar y^{\alpha_n},\bar y_\Gamma^{\alpha_n}, \bar w^{\alpha_n})$ {\em of the corresponding
state system} (\ref{eq:2.13})--(\ref{eq:2.16}), {\em the convergences}  (\ref{eq:3.18})--(\ref{eq:3.21}) {\em
hold as} $\,n\to\infty$. 

\vspace{2mm}
(ii) \,\,{\em Whenever sequences $\,\{\alpha_n\}\subset (0,1]\,$ and $\,\{(\bar y^{\alpha_n}, \bar y_\Gamma^{\alpha_n},  
\bar u_\Gamma^{\alpha_n})\}$ having the
properties described in\/} (i) {\em are given, then the following holds true: to any subsequence 
$\{n_k\}_{k\in\nz}$ of $\nz$ there are a {subsequence} $\,\{n_{k_\ell}\}_{\ell\in\nz}\,$ and some
$((\lambda,\lambda_\Gamma),(q,q_\Gamma))\in {\cal W}_0^*\times {\cal Z}\,$ such that }
\begin{itemize}
\item {\em the relations} (\ref{eq:4.35})--(\ref{eq:4.38}), (\ref{eq:4.42}), {\em and} (\ref{eq:4.44})
{\em hold (where the sequences are indexed by $\,n_{k_\ell}\,$ and the limits are taken {\pier as}
$\ell\to\infty$), and}  
\item {\em the variational inequality} (\ref{eq:4.39}) {\em and the adjoint equation} (\ref{eq:4.41})
{\em are satisfied.}
\end{itemize}

\vspace{3mm}
{\bf Remark~4.5:}\quad\,Unfortunately, we are not able to show that the limit pair $(q,q_\Gamma)$ solving
the adjoint problem associated with the optimal triple $(\bar{y},\bar{y}_\Gamma,\bar u_\Gamma)$ is uniquely determined. Therefore, it may well happen that the limiting pairs differ for different
subsequences. However, it follows from the variational inequality (\ref{eq:4.39}) 
that for any such limit pair {\Hassan $(q,q_\Gamma)$} it holds, with the orthogonal projection  ${\rm I\!P}_{\uad}$ 
onto $\uad$ with respect to the standard inner product in $H_\Gamma$, that for $\beta_5>0$ we have 
\begin{equation}
\label{eq:4.45}
\bar u_\Gamma={\rm I\!P}_{\uad}\left(-\beta_5^{-1}q_\Gamma\right)\,.
\end{equation}  
Standard arguments then yield that {\pier if the function $\bar{\bar u}_\Gamma \in L^2 (\Sigma) $ defined by  
\begin{equation}
\label{eq:4.46}
\bar{\bar u}_\Gamma (x,t)=\left\{
\begin{array}{ll}
\widetilde{u}_{2_\Gamma}(x,t)&\mbox{if } -\beta_5^{-1}q_\Gamma(x,t)>\widetilde{u}_{2_\Gamma}(x,t)\\[1mm]
\widetilde{u}_{1_\Gamma}(x,t)&\mbox{if } -\beta_5^{-1}q_\Gamma(x,t)<\widetilde{u}_{1_\Gamma}(x,t)\\[1mm]
-\beta_5^{-1}q_\Gamma(x,t)&\mbox{otherwise}
\end{array}
\right. \quad\hbox{for a.\,a. } (x,t)\in\Sigma \, , 
\end{equation}
belongs to $\uad$ (i.e., its time derivative actually exists and satisfies the bound prescribed in (\ref{eq:1.8})), then 
$\bar{\bar u}_\Gamma = {\bar u}_\Gamma $  and ${\bar u}_\Gamma $ turns out to be a pointwise projection.}
%%%%%%%%%%%%%%%%%%%%%%%%%%%%%%%%%%%%%%%%%%%%%%%%%%%%%%%%
 
\end{document}